\documentclass[11pt,a4paper,leqno]{article}
\usepackage{amscd,amsmath,amsthm}
\usepackage{pb-diagram}
\newtheorem{theorem}{Theorem}[section]
\newtheorem{lemma}[theorem]{Lemma}
\newtheorem{corollary}[theorem]{Corollary}
\newtheorem{proposition}[theorem]{Proposition}

\newtheorem{fact}[theorem]{Fact}
\newtheorem{definition}{Definition}[section]
\newtheorem{example}{Example}[section]
 \input amssym.def
\input amssym.tex

\def\bc{\begin{center}}
\def\ec{\end{center}}

\begin{document}
\title{Characterization Of Left Artinian Algebras Through \\Pseudo Path Algebras}
\author{ Fang Li}

 \date{}
\maketitle
\begin{abstract}

 In this paper, using pseudo path algebras, we generalize Gabriel's Theorem
on elementary
 algebras to left Artinian algebras over a field $k$ when it is splitting over its radical,
  in particular, when the dimension of the quotient algebra decided by the $n$'th Hochschild
  cohomology is less than $2$ (for example,  $k$ is finite or char$k=0$). Using generalized path algebras, the generalized
Gabriel's Theorem is given for finite dimensional algebras with
2-nilpotent radicals which is splitting over its radical. As a
tool, the so-called pseudo path algebras are introduced as a new
generalization of path algebras, which can cover generalized path
algebras (see Fact 2.5).

 The main result is that (i) for a left Artinian $k$-algebra $A$ and $r=r(A)$ the radical of $A$, when
the quotient algebra $A/r$ can be lifted, it holds that
 $A\cong PSE_k(\Delta,\mathcal{A},\rho)$ with
$J^{s}\subset\langle\rho\rangle\subset J$ for some $s$ (Theorem
3.2); (ii) for a finite dimensional $k$-algebra $A$ with $r=r(A)$
2-nilpotent radical, when the quotient algebra $A/r$ can be
lifted, it holds that
  $A\cong k(\Delta,\mathcal{A},\rho)$ with
$\widetilde J^{2}\subset\langle\rho\rangle\subset\widetilde
J^{2}+\widetilde J\cap$ \textrm{Ker}$\widetilde{\varphi}$ (Theorem
4.3), where $\Delta$ is the quiver of $A$ and $\rho$ is a set of
relations.

 Meantime,  the
uniqueness of the quivers and generalized path algebra/pseudo path
algebras satisfying the isomorphism relations
   is obtained in the case when the ideals generated by the relations are admissible (see Theorem 3.5 and 4.4).
\end{abstract}
\textbf{2000 Mathematics Subject Classifications:} 16G10

\section{Introduction}

In this paper, $k$ will always denote a field, and all modules are
unital. An algebra is said to be  {\em left Artinian} if it
satisfies the descending chain condition on left ideals.

It is well-known that for a finite dimensional algebra $A$ over
 an algebraically closed field $k$ and
 the nilpotent radical $N=J(A)$, the quotient algebra $A/N$ is
 semisimple, that is, there are uniquely positive integers
 $n_{1}\leq n_{2}\leq\cdot\cdot\cdot\leq n_{r}$ such that $A/N\cong
 M_{n_{1}}(k)\oplus\cdot\cdot\cdot\oplus M_{n_{r}}(k)$ where
 $M_{n_{i}}(k)$ denote the algebras of $n_{i}\times n_{i}$
 matrices with entries in $k$, which are trivially $k$-simple algebras.  In the special case that $A$ is
 an  {\em elementary} algebra$^{\cite{ARS}}$,  every $n_{i}=1$,
 that is,
 $M_{n_{i}}\cong k$, or say, $A/N$ is a direct sum of some $k$ as $k$-algebras,
 writing $A/N=\coprod_{r}(k)$.

Obviously, every finite dimensional path algebra is elementary.
Conversely, by the famous Gabriel's Theorem$^{\cite{ARS}}$, for
each elementary algebra $\Lambda$, one can construct the
correspondent quiver $\Gamma(\Lambda)$ of $\Lambda$
 such that $\Lambda$ is isomorphic to a quotient algebra of the path
 algebra $k\Gamma(\Lambda)$. On the other hand,  the
 module category of any algebra $A$ always is Morita-equivalent to that of some elementary
 algebra$^{\cite{DK}}$.
Therefore, in the view point of representation theory, it should
 be enough to consider representations of elementary
algebras, or say equivalently, quotient algebras of path algebras.
In particular, it has provided the description of finitely
generated modules over some given algebras (for instance
\cite{ARS}\cite{R}).

However, in the view point of structures of algebras, finite
dimensional algebras can not be replaced with elementary algebras,
for example,   if one should consider how to give the possible
classification for finite dimensional algebras.

In this reason, Shao-xue Liu, one of the authors of \cite{CL},
raised an interesting problem, that is, how to find a
generalization of path algebras so as to obtain a generalized type
of the Gabriel Theorem for arbitrary finite dimensional algebras
which would admit this algebra to be isomorphic to a quotient
algebra of such a generalized path algebra. The first step on this
idea was finished in \cite{CL} where the valid concept of
generalized path algebras was introduced (see Section 2). But, the
further research could not be found.

In this paper, we hope to answer the Liu's problem affirmatively
by using of pseudo path algebras and generalized path algebras
under the meaning of \cite{CL}.

Some preparation is given in Section 2. In fact, we find that
generalized path algebras is not enough to characterize algebras
unless they are finite-dimensional with $2$-nilpotent radicals. In
this reason, the so-called pseudo path algebras are introduced as
a new generalization of path algebras, which can cover generalized
path algebras (see Fact 2.5).
 In Section 3, using of pseudo path algebras, we generalize the Gabriel's Theorem
on elementary
 algebras to left Artinian algebras over a field $k$ in the case that the quotient algebra is lifted,
  in particular, when the Dimension of the quotient algebra decided by the $n$'th Hochschild
  cohomology is less than $2$ (for example,  $k$ is finite or char$k=0$).
On the other hand, in Section 4, relying on generalized path
algebras, the Generalized Gabriel's Theorem is given for finite
dimensional algebras with 2-nilpotent radicals in the case that
the quotient algebra is lifted.
 In all cases we discuss,  the uniqueness of such quivers $\Delta$ and generalized path algebras/pseudo path algebras
   is taken on in the case that the ideals generated by the relations are admissible (see Theorem 3.5 and 4.4).

Under some certain conditions, the generalized Gabriel's Theorems
are not depended on the ground field. This fact gives the
possibility to provide an approach to modular representations of
algebras and groups.

Note that when
 $A\cong k(\Delta,\mathcal{A})/\langle\rho\rangle$ or $A\cong PSE_k(\Delta,\mathcal{A})/\langle\rho\rangle$,  the
 structure of $A$ is decided by the generated ideal
  $\langle\rho\rangle$ of the set $\rho$ of some relations. From
  this, one can try to classify those associative algebras satisfying
  the theorems, including many important kinds of algebras.
   All these above are the works we are processing after this
   paper, which will present the further significance.

\section{On Generalized Path Algebras And Pseudo Path Algebras}

In this section, we firstly introduce the definitions of
generalized path algebra$^{\cite{CL}}$ and pseudo path algebra and
then discuss their properties and relationship.

A {\em quiver} $\Delta$ is given by two sets $\Delta_{0}$ and
$\Delta_{1}$ together with two maps $s, e$:
$\Delta_{1}\rightarrow\Delta_{0}$. The elements of $\Delta_{0}$
are called {\em vertices}, while the elements of $\Delta_{1}$ are
called {\em arrows}. For an arrow $\alpha\in\Delta_{1}$, the
vertex $s(\alpha)$ is the {\em start vertex} of $\alpha$ and the
vertex $e(\alpha)$ is the {\em end vertex } of $\alpha$, and we
draw $s(\alpha)\stackrel{\alpha}{\rightarrow}e(\alpha)$. A {\em
path} $p$ in $\Delta$ is $(a|\alpha_{1}\cdots\alpha_{n}|b)$, where
$\alpha_{i}\in\Delta_{1}$, for $i=1,...,n$, and $s(\alpha_{1})=a$,
$e(\alpha_{i})=s(\alpha_{i+1})$ for $i=1,...,n+1$, and
$e(\alpha_{n})=b$.
 $s(\alpha_{1})$ and $e(\alpha_{n})$ are also called respectively  the start vertex and the end vertex of
 $p$. Write $s(p)=s(\alpha_{1})$ and $e(p)=e(\alpha_{n})$.
  The {\em length} of a path is the number of arrows in it. To each
arrow $\alpha$, one can assign an edge $\overline{\alpha}$ where
the orientation is forgotten. A {\em walk} between two vertices
$a$ and $b$ is given by
$(a|\overline{\alpha_{1}}\cdot\cdot\cdot\overline{\alpha_{n}}|b)$,
where $a\in\{s(\alpha_{1}),e(\alpha_{1})\}$,
$b\in\{s(\alpha_{n}),e(\alpha_{n})\}$, and for each $i=1, . . .
,n-1$,
$\{s(\alpha_{i}),e(\alpha_{i})\}\bigcap\{s(\alpha_{i+1}),e(\alpha_{i+1})\}\not=\emptyset$.
A quiver is said to be {\em connected} if there exists a walk
between any two vertices $a$ and $b$.

In this paper, we will always assume the quiver $\Delta$ is
finite, i.e. the number of vertices $|\Delta_0|<\infty$ and
$|\Delta_1|<\infty$.

\begin{definition}\label{def1.1}
For two algebras $A$ and $B$,  the {\em rank} of a finitely
generated $A$-$B$-bimodule $M$ is defined as the least cardinal
number of the sets of generators. In particular, if  $M=0$, it is
said to be with rank $0$ as a finitely generated $A$-$B$-bimodule.
\end{definition}

Clearly, for any finitely generated $A$-$B$-bimodule, such rank
always exists uniquely.
\\

{\bf Generalized Path Algebra And Tensor Algebra}

Let $\Delta=(\Delta_{0},\Delta_{1})$ be a quiver and
$\mathcal{A}=\{A_{i}: i\in\Delta_{0}\}$ be a family of $k$-algebra
$A_{i}$ with identity $e_{i}$, indexed by the vertices of
$\Delta$. The elements $a_{i}$ of $\bigcup_{i\in\Delta_{0}}A_{i}$
are called  the {\em $\mathcal{A}$-paths of length zero}, whose
start vertex $s(a_{i})$ and the end vertex $e(a_{i})$ are   both
$i$. For each $n\geq 1$, an {\em $\mathcal{A}$-path $P$ of length
$n$} is given by $a_{1}\beta_{1}a_{2}\beta_{2}\cdot\cdot\cdot
a_{n}\beta_{n}a_{n+1}$, where
$(s(\beta_{1})|\beta_{1}\cdot\cdot\cdot\beta_{n}|e(\beta_{n}))$ is
a path in $\Delta$ of length $n$, for each $i=1,...,n$, $a_{i}\in
A_{s(\beta_{i})}$ and $a_{n+1}\in A_{e(\beta_{n})}$.
  $s(\beta_{1})$ and $e(\beta_{n})$ are also called respectively  the start vertex and the end vertex of
 $P$. Write $s(P)=s(\alpha_{1})$ and $e(P)=e(\alpha_{n})$.
 Now, consider the quotient $R$ of
the $k$-linear space with basis the set of all $\mathcal{A}$-paths
of $\Delta$ by the subspace generated by all the elements of the
form
\[
a_{1}\beta_{1}\cdot\cdot\cdot\beta_{j-1}(a_{j}^{1}+\cdot\cdot\cdot+a_{j}^{m})\beta_{j}a_{j+1}\cdot\cdot\cdot\beta_{n}a_{n+1}
-\sum_{l=1}^{m}a_{1}\beta_{1}\cdot\cdot\cdot\beta_{j-1}a_{j}^{l}\beta_{j}a_{j+1}\cdot\cdot\cdot\beta_{n}a_{n+1}
\]
where
$(s(\beta_{1})|\beta_{1}\cdot\cdot\cdot\beta_{n}|e(\beta_{n}))$ is
a path in $\Delta$ of length $n$, for each $i=1,...,n$, $a_{i}\in
A_{s(\beta_{i})}$, $a_{n+1}\in A_{e(\beta_{n})}$ and $a_{j}^{l}\in
A_{s(\beta_{j})}$ for $l=1,...,m$. In $R$,  given two elements
 $[a_{1}\beta_{1}a_{2}\beta_{2}\cdot\cdot\cdot a_{n}\beta_{n}a_{n+1}]$ and
 $[b_{1}\gamma_{1}b_{2}\gamma_{2}\cdot\cdot\cdot b_{n}\gamma_{n}b_{n+1}]$ in
 $R$, define the multiplication as follows:
\\
\\
$ [a_{1}\beta_{1}a_{2}\beta_{2}\cdot\cdot\cdot
a_{n}\beta_{n}a_{n+1}]\cdot[b_{1}\gamma_{1}b_{2}\gamma_{2}\cdot\cdot\cdot
b_{n}\gamma_{n}b_{n+1}]\\
\\
=\left\{\begin{array}{ll}
 [a_{1}\beta_{1}a_{2}\beta_{2}\cdot\cdot\cdot
a_{n}\beta_{n}(a_{n+1}b_{1})\gamma_{1}b_{2}\gamma_{2}\cdot\cdot\cdot
b_{n}\gamma_{n}b_{n+1}],
  &  \mbox{if $a_{n+1}, b_{1}\in A_{i}$ for same $i$}\\
0,  &  \mbox{otherwise}
\end{array}
\right. $ \\
\\
It is easy to check that the above multiplication is well-defined
 and makes $R$ to become a $k$-algebra. This algebra $R$
defined above is called the {\em $\mathcal{A}$-path algebra} of
$\Delta$. Denote it by $R=k(\Delta,\mathcal{A})$. Clearly, $R$ is
an $A$-bimodule, where $A=\oplus_{i\in\Delta_0}A_{i_0}$. All such
algebras are said to be {\em generalized path algebras}.

Remark that (i) $R=k(\Delta,\mathcal{A})$ has identity if and only
if $\Delta_{0}$ is finite; (ii) Any path
$(s(\beta_{1})|\beta_{1}\cdot\cdot\cdot\beta_{n}|e(\beta_{n}))$ in
$\Delta$ can be considered as an $\mathcal{A}$-path with
$a_{i}=e_{i}$. Hence the usual path algebra $k\Delta$ can be
embedded into the $\mathcal{A}$-path algebra
$k(\Delta,\mathcal{A})$. Or say, if $A_{i}=k$ for each
$i\in\Delta_{0}$, then $k(\Delta,\mathcal{A})=k\Delta$; (iii) For
$R=k(\Delta,\mathcal{A})$, $dim_{k}R<\infty$ if and only if
$dim_{k}A_{i}<\infty$ for each $i\in\Delta_{0}$ and $\Delta$ is a
finite quiver without oriented cycles.

Associated with the pair $(A, _{A}M_{A})$ for a $k$-algebra $A$
and an $A$-bimodule $M$, we write the $n$-fold $A$-tensor product
$M\otimes_{A}M\otimes\cdot\cdot\cdot\otimes_{A}M$ as $M^{n}$, then
$T(A,M)=A\oplus M\oplus M^{2}\oplus\cdot\cdot\cdot\oplus
M^{n}\oplus\cdot\cdot\cdot$ as an abelian group. Writing
$M^{0}=A$, then $T(A,M)$ becomes a $k$-algebra with multiplication
induced by the natural $A$-bilinear maps $M^{i}\times
M^{j}\rightarrow M^{i+j}$ for $i\geq 0$ and $j\geq 0$.  $T(A,M)$
is called the {\em tensor algebra} of $M$ over $A$.

Now, we define a special class of tensor algebras so as to
characterize generalized path algebras. An {\em
$\mathcal{A}$-path-type tensor algebra} is defined to be the
tensor algebra $T(A,M)$ satisfying that (i)
$A=\bigoplus_{i\in\Delta_{0}}A_{i}$ for a family of $k$-algebras
$\mathcal{A}=\{A_{i}: i\in\Delta_{0}\}$, (ii) $M=\bigoplus_{i,j\in
I}$$ _{i}M_{j}$ for finitely generated $A_{i}$-$A_{j}$-bimodules
$_{i}M_{j}$ for all $i$ and $j$ in $I$ and $A_{k}\cdot_{i}M_{j}=0$
if $k\not=i$ and
 $_{i}M_{j}\cdot A_{k}=0$ if $k\not=j$.
 A {\em free
$\mathcal{A}$-path-type tensor algebra} is the
$\mathcal{A}$-path-type tensor algebra $T(A,M)$ whose each
finitely generated $A_{i}$-$A_{j}$-bimodule $_{i}M_{j}$ for $i$
and $j$ in $I$ is a free bimodule with a basis and the rank of
this basis is equal to the rank of $_{i}M_{j}$ as a finitely
generated $A_{i}$-$A_{j}$-bimodule.

 $\mathcal{A}$-path-type tensor algebras and generated path
 algebras can be constructed each other as follows.

For an $\mathcal{A}$-path algebra $k(\Delta,\mathcal{A})$, let
$A=\bigoplus_{i\in\Delta_{0}}A_{i}$. For any $i$ and $j$, let
$_{i}M^{F}_{j}$ be the free $A_{i}$-$A_{j}$-bimodule with basis
given by the arrows from $i$ to $j$. It is easy to see that the
number of free generators in the basis is the rank of
$_{i}M^{F}_{j}$ as a finitely generated bimodule. Define
$A_{k}\cdot_{i}M^{F}_{j}=0$ if $k\not=i$ and
 $_{i}M^{F}_{j}\cdot A_{k}=0$ if $k\not=j$. Let $M^{F}=\bigoplus_{i\rightarrow j}$$ _{i}M^{F}_{j}$, which is
 clearly an $A$-bimodule. Then we  get uniquely the free $\mathcal{A}$-path-type tensor algebras $T(A,M^{F})$.

Conversely, assume that $T(A,M)$ is an $\mathcal{A}$-path-type
tensor algebra with a family of $k$-algebras $\mathcal{A}=\{A_{i}:
i\in I\}$ and finitely generated $A_{i}$-$A_{j}$-bimodules
$_{i}M_{j}$ for $i,j\in I$ such that $A=\bigoplus_{i\in I}A_{i}$,
$M=\bigoplus_{i,j\in I}$$ _{i}M_{j}$, $A_{k}\cdot_{i}M_{j}=0$ if
$k\not=i$ and
 $_{i}M_{j}\cdot A_{k}=0$ if $k\not=j$. Trivially, $_{i}M_{j}=A_{i}MA_{j}$.
 Let the rank of $_{i}M_{j}$ be $r_{ij}$.
Now we can associate with $T(A,M)$ a quiver
$\Delta=(\Delta_{0},\Delta_{1})$ and its generalized path algebra
$R=k(\Delta,\mathcal{A})$ in the following way. Let $\Delta_{0}=I$
as the set of vertices. For $i,j\in I$, let the number of arrows
from $i$ to $j$ in $\Delta$ be the rank $r_{ij}$ of the finitely
generated $A_{i}$-$A_{j}$-bimodules $_{i}M_{j}$. Obviously, if
$_{i}M_{j}=0$, then there are no arrows from $i$ to $j$. Thus, we
get a quiver $\Delta=(\Delta_{0},\Delta_{1})$ which is called {\em
the quiver of $T(A,M)$}, and its $\mathcal{A}$-path algebra
$R=k(\Delta,\mathcal{A})$ which is called {\em the corresponding
$\mathcal{A}$-path algebra of $T(A,M)$}.

One can find two non-isomorphic finitely generated bimodules
 which possesses the same ranks, therefore
there are two $\mathcal{A}$-path-type tensor algebras $T(A,M_{1})$
and $T(A,M_{2})$, with non-isomorphic bimodules $M_{1}$ and
$M_{2}$, but their induced quivers and $\mathcal{A}$-path algebras
are the same in the above way.

From the above discussion, every $\mathcal{A}$-path-type tensor
algebra $T(A,M)$ can be used to construct its corresponding
$\mathcal{A}$-path algebra $k(\Delta,\mathcal{A})$; but, from this
$\mathcal{A}$-path algebra $k(\Delta,\mathcal{A})$, we can get
uniquely the free $\mathcal{A}$-path-type tensor algebra
$T(A,M^{F})$. Thus, we have:

\begin{lemma}
Every $\mathcal{A}$-path-type tensor algebra $T(A,M)$ can be used
to construct uniquely the free $\mathcal{A}$-path-type tensor
algebra $T(A,M^{F})$. There is a surjective $k$-algebra morphism
$\pi$: $T(A,M^{F})\rightarrow T(A,M)$ such that
$\pi(_{i}M^{F}_{j})=_{i}M_{j}$ for any $i,j\in I$.
\end{lemma}
{\em Proof}: It is sufficient to prove the second conclusion. For
$T(A,M)$, let the rank of $_{i}M_{j}$ be $r_{ij}$. Thus, for the
corresponding $\mathcal{A}$-path algebra $k(\Delta,\mathcal{A})$,
the number of the arrows from $i$ to $j$ is $r_{ij}$, and then, in
$T(A,M^{F})$,  the rank of the free generators of $_{i}M^{F}_{j}$
given by the arrows is also $r_{ij}$. Define $\pi$:
$T(A,M^{F})\rightarrow T(A,M)$ by giving a bijection between the
set of the free generators of $_{i}M^{F}_{j}$ and the set of the
chosen generators of $_{i}M_{j}$ with number of the rank. Then
$\pi$ can be expanded to become into a surjective $k$-algebra
morphism with $\pi(_{i}M^{F}_{j})=_{i}M_{j}$ for any $i,j\in I$.
$\;\;\;\;\square$

Next, we will show in the following Proposition 2.9 that every
$\mathcal{A}$-path-type tensor algebra is a homomorphic image of
its corresponding $\mathcal{A}$-path algebra.

The following criterion (that is, Lemma III.1.2 in \cite{ARS})  is
useful for constructing algebra morphisms from tensor algebras to
other algebras.
\begin{lemma}
Let $A$ be a $k$-algebra and $M$ an $A$-bimodule. Let $\Lambda$ be
a $k$-algebra and $f: A\oplus M\rightarrow\Lambda$ a map such that
the following two conditions are satisfied:

(i) $f|_{A}: A\rightarrow\Lambda$ is an algebra morphism;

 (ii) Viewing $f(M)$ as an $A$-bimodule via $f|_{A}:
 A\rightarrow\Lambda$ then $f|_{M}: M\rightarrow f(M)\subset\Lambda$ is an
 $A$-bimodule map.\\
 Then there is a unique algebra morphism $\widetilde{f}:
 T(A,M)\rightarrow\Lambda$ such that  $\widetilde{f}|_{A\oplus
 M}=f$ and generally, $\widetilde{f}(\sum_{n=0}^{\infty}m^{n}_{1}\otimes\cdot\cdot\cdot\otimes m^{n}_{n})
 =\sum_{n=0}^{\infty}f(m^{n}_{1})\cdot\cdot\cdot f(m^{n}_{n})$ for $m^{n}_{1}\otimes\cdot\cdot\cdot\otimes m^{n}_{n}\in M^{n}$.
\end{lemma}

Note that it is enough for the proof of (ii) in \cite{ARS} under
the condition that $f(M)$ is an $A$-bimodule via $f|_{A}:
 A\rightarrow\Lambda$.

Clearly, all $\mathcal{A}$-paths of length zero, that is, the
elements of $\bigcup_{i\in\Delta_{0}}A_{i}$ can generate a
subalgebra of $k(\Delta,\mathcal{A})$, which is denoted by
$k(\Delta_{0},\mathcal{A})$. And, denote by
$k(\Delta_{1},\mathcal{A})$ the $k$-linear space consisting of all
$\mathcal{A}$-paths of length $1$ and $J$ the ideal in a
$\mathcal{A}$-path algebra $k(\Delta,\mathcal{A})$ generated by
all elements in $k(\Delta_{1},\mathcal{A})$. It is easy to see
that $k(\Delta_{1},\mathcal{A})$ is an $A$-subbimodule of
$k(\Delta,\mathcal{A})$.
\\

{\bf Pseudo Path Algebra And Pseudo Tensor Algebra}

Let $\Delta=(\Delta_{0},\Delta_{1})$ be a quiver and
$\mathcal{A}=\{A_{i}: i\in\Delta_{0}\}$ be a family of $k$-algebra
$A_{i}$ with identity $e_{i}$, indexed by the vertices of
$\Delta$. The elements $a_{i}$ of $\bigcup_{i\in\Delta_{0}}A_{i}$
are called  the {\em $\mathcal{A}$-pseudo-paths of length zero},
whose start vertex $s(a_{i})$ and the end vertex $e(a_{i})$  both
are $i$. For each $n\geq 1$, a {\em pure $\mathcal{A}$-pseudo-path
$P$ of length $n$} is given by $a_{1}\beta_{1}b_{1}\cdot
a_{2}\beta_{2}b_{2}\cdot...\cdot a_{n}\beta_{n}b_{n}$, where
$(s(\beta_{1})|\beta_{1}\cdot\cdot\cdot\beta_{n}|e(\beta_{n}))$ is
a path in $\Delta$ of length $n$, for each $i=1,...,n$,
$b_{i-1}\in A_{e(\beta_{i-1})}, a_{i}\in A_{s(\beta_{i})}$ with
$s(\beta_{i})=e(\beta_{i-1})$.
 $s(\beta_{1})$ and $e(\beta_{n})$ are also called respectively  the {\em start vertex} and the {\em end vertex}
  of $P$. Write $s(P)=s(\beta_{1})$ and $e(P)=e(\beta_{n})$. A  {\em
general $\mathcal{A}$-pseudo-path $Q$ of length $n$} is given in
the form

$\alpha_{1}\cdot c_{1}\cdot \alpha_{2}\cdot c_{2}\cdot...\cdot
c_{k}\cdot\alpha_{k}$\\or

$c_{0}\cdot\alpha_{1}\cdot c_{1}\cdot \alpha_{2}\cdot
c_{2}\cdot...\cdot c_{k}\cdot\alpha_{k}$\\or

 $\alpha_{1}\cdot
c_{1}\cdot \alpha_{2}\cdot c_{2}\cdot...\cdot
c_{k}\cdot\alpha_{k}\cdot c_{k+1}$\\
 or

 $c_{0}\cdot\alpha_{1}\cdot
c_{1}\cdot \alpha_{2}\cdot c_{2}\cdot...\cdot
c_{k}\cdot\alpha_{k}\cdot c_{k+1}$ \\
where $\alpha_{i}$ is an pure $\mathcal{A}$-pseudo-path of length
$n_{i}$ and $\sum_{i=1}^{k}n_{i}=n$, and the start vertex of
$\alpha_{i+1}$ is just the end vertex of $\alpha_{i}$, that is,
$e(\alpha_{i})=s(\alpha_{i+1})$ and $c_{i}\in A_{e(\alpha_{i})}$.

Let $V$ be the $k$-linear space with basis the set of all general
$\mathcal{A}$-paths of $\Delta$.

  Consider the quotient $R$ of
the $k$-linear space $V$ by the subspace generated by all the
elements of the form
\begin{equation}
a_{1}\beta_{1}b_{1}\cdot...\cdot
a_{j}\beta_{j}(b_{j}^{1}+...+b_{j}^{m})\cdot\gamma
-\sum_{l=1}^{m}a_{1}\beta_{1}b_{1}\cdot...\cdot
a_{j}\beta_{j}b_{j}^{l}\cdot\gamma
\end{equation}
\begin{equation}
\alpha\cdot(a_{1}^{1}+...+a_{1}^{m}) \beta_{1}b_{1}\cdot...\cdot
a_{n}\beta_{n}b_{n} -\sum_{l=1}^{m}\alpha\cdot a_{1}^{l}
\beta_{1}b_{1}\cdot...\cdot a_{n}\beta_{n}b_{n}
\end{equation}
\begin{equation}
(ab)\cdot c\beta d-a\cdot(b\cdot c\beta d),\>\>\> a\beta
b\cdot(cd)-(a\beta b\cdot c)\cdot d
\end{equation}
\begin{equation}
a\beta b\cdot 1-a\beta b,\>\>\> 1\cdot a\beta b-a\beta b
\end{equation}\\
where
$a,b,c,d,b_{j}^{l},a_{1}^{l}\in\bigcup_{i\in\Delta_{0}}A_{i}$ and
$1$ is the identity of $A=\oplus_{i\in\Delta_0}A_i$.

 In $R$, define the following multiplication. Given two elements
 $[a_{1}\beta_{1}b_{1}\cdot a_{2}\beta_{2}b_{2}\cdot...\cdot a_{n}\beta_{n}b_{n}]$ and
$[c_{1}\gamma_{1}d_{1}\cdot c_{2}\gamma_{2}d_{2}\cdot...\cdot
c_{n}\gamma_{m}d_{m}]$ in which at least one is of length $n\geq
1$, define\\
\\
$[a_{1}\beta_{1}b_{1}\cdot a_{2}\beta_{2}b_{2}\cdot...\cdot
a_{n}\beta_{n}b_{n}]\cdot[c_{1}\gamma_{1}d_{1}\cdot
c_{2}\gamma_{2}d_{2}\cdot...\cdot c_{n}\gamma_{m}d_{m}]\\
\\
=\left\{\begin{array}{ll} [a_{1}\beta_{1}b_{1}\cdot
a_{2}\beta_{2}b_{2}\cdot...\cdot a_{n}\beta_{n}b_{n}\cdot
c_{1}\gamma_{1}d_{1}\cdot c_{2}\gamma_{2}d_{2}\cdot...\cdot
c_{n}\gamma_{m}d_{m}] ,
  &  \mbox{if $b_{n}, c_{1}\in A_{i}$ for same $i$}\\
0,  &  \mbox{otherwise.}
\end{array}
\right. $ \\
\\
Given two elements $a, b$ of length zero, that is,
$a,b\in\bigcup_{i\in\Delta_{0}}A_{i}$, define \\
$a\cdot b=\left\{\begin{array}{ll} ab, & \mbox{if $a, b$ are in
the same
$A_{i}$}\\
0,  &  \mbox{otherwise}
\end{array}
\right. $ $\;\;\;$ where $ab$ means the product of $a, b$ in
$A_i$.

It is easy to check that the above multiplication in $R$ is
well-defined and makes $R$ to become a $k$-algebra. This algebra
$R$ defined above is called the {\em $\mathcal{A}$-pseudo path
algebra} of $\Delta$. Denote it by
$R=PSE_{k}(\Delta,\mathcal{A})$. Clearly, $R$ is an $A$-bimodule.

Remark that (i) $R=PSE_{k}(\Delta,\mathcal{A})$ has identity if
and only if $\Delta_{0}$ is finite; (ii) Any path
$(s(\beta_{1})|\beta_{1}\cdot\cdot\cdot\beta_{n}|e(\beta_{n}))$ in
$\Delta$ can be considered as an $\mathcal{A}$-path with
$a_{i}=e_{i}$ the identity of $A_{i}$. Hence the usual path
algebra $k\Delta$ can be embedded into the $\mathcal{A}$-pseudo
path algebra $PSE_k(\Delta,\mathcal{A})$. Or say, if $A_{i}=k$ for
each $i\in\Delta_{0}$, then $PSE_k(\Delta,\mathcal{A})=k\Delta$;
(iii) For $R=PSE_k(\Delta,\mathcal{A})$, $dim_{k}R<\infty$ if and
only if $dim_{k}A_{i}<\infty$ for each $i\in\Delta_{0}$ and
$\Delta$ is a finite quiver without oriented cycles.

Associated with the pair $(A, _{A}M_{A})$ for a $k$-algebra $A$
and an $A$-bimodule $M$, we write the $n$-fold $k$-tensor product
$M\otimes_{k}M\otimes\cdot\cdot\cdot\otimes_{k}M$ as $M^{n}$;
$\sum_{M_{1},M_{2},\cdot\cdot\cdot,M_{n}}M_{1}\otimes_{k}M_{2}\otimes_{k}\cdot\cdot\cdot\otimes_{k}M_{n}$
as $M(n)$ where $M_{i}$ is either $M$ (at least there exists one)
or $A$ but no two $A$'s are neighbouring, then define
$\mathcal{PT}(A,M)=A\oplus M(1)\oplus
M(2)\oplus\cdot\cdot\cdot\oplus M(n)\oplus\cdot\cdot\cdot$ as an
abelian group. Denote by $M(n,l)$ the sum of these items
$M_{1}\otimes_{k}M_{2}\otimes_{k}\cdot\cdot\cdot\otimes_{k}M_{n}$
of $M(n)$ in which there are $l$ $M_{i}$'s equal to $M$. Clearly,
$(n-1)/2\leq l\leq n$ and $M(n)=\sum_{(n-1)/2\leq l\leq n}M(n,l)$.
Writing $M^{0}=A$, then $\mathcal{PT}(A,M)$ becomes a $k$-algebra
with multiplication induced by the natural $k$-bilinear maps:

 $M^{i}\times M^{j}\rightarrow M^{i+j}$ $\>\>\>\>$for $i\geq 1$,$\>\>$ $j\geq 1$;

 $M^{i}\times A\rightarrow M^{i}\otimes_{k}A$ $\>\>\>\>$for $i\geq 1$;

 $A\times M^{j}\rightarrow A\otimes_{k}M^{j}$ $\>\>\>\>$for $j\geq
 1$\\
and the natural $A$-bilinear map:

 $A\times A\rightarrow A\otimes_{A}A=A$.\\
 Note that the associative law of $\mathcal{PT}(A,M)$ is from the
 fact $(A\otimes_{A}A)\otimes_{k}M\cong
 A\otimes_{A}(A\otimes_{k}M)$.
  We call $\mathcal{PT}(A,M)$ a {\em pseudo tensor algebra}.

Now, we define a special class of pseudo tensor algebras so as to
characterize pseudo path algebras. An {\em $\mathcal{A}$-path-type
pseudo tensor algebra} is defined to be the pseudo tensor algebra
$\mathcal{PT}(A,M)$ satisfying that (i)
$A=\bigoplus_{i\in\Delta_{0}}A_{i}$ for a family of $k$-algebras
$\mathcal{A}=\{A_{i}: i\in\Delta_{0}\}$, (ii) $M=\bigoplus_{i,j\in
I}$$ _{i}M_{j}$ for finitely generated $A_{i}$-$A_{j}$-bimodules
$_{i}M_{j}$ for all $i$ and $j$ in $I$ and $A_{k}\cdot_{i}M_{j}=0$
if $k\not=i$ and
 $_{i}M_{j}\cdot A_{k}=0$ if $k\not=j$.
 A {\em free
$\mathcal{A}$-path-type pseudo tensor algebra} is the
$\mathcal{A}$-path-type pseudo tensor algebra $\mathcal{PT}(A,M)$
whose each finitely generated $A_{i}$-$A_{j}$-bimodule $_{i}M_{j}$
for $i$ and $j$ in $I$ is a free bimodule with a basis and the
rank of this basis is equal to the rank of $_{i}M_{j}$ as a
finitely generated $A_i$-$A_j$-bimodule.

 $\mathcal{A}$-path-type pseudo tensor algebras and pseudo path
 algebras can be constructed each other as follows.

Given an $\mathcal{A}$-pseudo path algebra
$PSE_k(\Delta,\mathcal{A})$, let
$A=\bigoplus_{i\in\Delta_{0}}A_{i}$. For any $i$ and $j$, let
$_{i}M^{F}_{j}$ be the free $A_{i}$-$A_{j}$-bimodule with basis
given by the arrows from $i$ to $j$. It is easy to see that the
number of free generator in the basis is the rank of
$_{i}M^{F}_{j}$ as a finitely generated bimodule. Define
$A_{k}\cdot_{i}M^{F}_{j}=0$ if $k\not=i$ and
 $_{i}M^{F}_{j}\cdot A_{k}=0$ if $k\not=j$. Let $M^{F}=\bigoplus_{i\rightarrow j}$$ _{i}M^{F}_{j}$, which is
 clearly an $A$-bimodule. Then we  get uniquely the free $\mathcal{A}$-path-type pseudo
  tensor algebras $\mathcal{PT}(A,M^{F})$.

Conversely, assume that $\mathcal{PT}(A,M)$ is an
$\mathcal{A}$-path-type pseudo tensor algebra with a family of
$k$-algebras $\mathcal{A}=\{A_{i}: i\in I\}$ and finitely
generated $A_{i}$-$A_{j}$-bimodules $_{i}M_{j}$ for all $i$ and
$j$ in $I$ such that $A=\bigoplus_{i\in I}A_{i}$ and
$M=\bigoplus_{i,j\in I}$$ _{i}M_{j}$, $A_{k}\cdot_{i}M_{j}=0$ if
$k\not=i$ and
 $_{i}M_{j}\cdot A_{k}=0$ if $k\not=j$. Trivially, $_{i}M_{j}=A_{i}MA_{j}$.
 Let the rank of $_{i}M_{j}$ be $r_{ij}$.
Now we can associate with $\mathcal{PT}(A,M)$ a quiver
$\Delta=(\Delta_{0},\Delta_{1})$ and its pseudo path algebra
$R=PSE_k(\Delta,\mathcal{A})$ in the following way. Let
$\Delta_{0}=I$ as the set of vertices. For $i,j\in I$, let the
number of arrows from $i$ to $j$ in $\Delta$ be the rank $r_{ij}$
of the finitely generated $A_{i}$-$A_{j}$-bimodules $_{i}M_{j}$.
Obviously, if $_{i}M_{j}=0$, then there are no arrows from $i$ to
$j$. Thus, we get a quiver $\Delta=(\Delta_{0},\Delta_{1})$ which
is called {\em the quiver of $\mathcal{PT}(A,M)$}, and its
$\mathcal{A}$-pseudo path algebra $R=PSE_k(\Delta,\mathcal{A})$
which is called {\em the corresponding $\mathcal{A}$-pseudo path
algebra of $\mathcal{PT}(A,M)$}.

One can find two non-isomorphic finitely generated bimodules which
possesses the same ranks,  therefore there are two
 $\mathcal{A}$-path-type pseudo tensor algebras
 $\mathcal{PT}(A,M_{1})$ and $\mathcal{PT}(A,M_{2})$, with
 non-isomorphic $M_{1}$ and $M_{2}$, but their induced
 quivers and  $\mathcal{A}$-pseudo path algebras are the same in the way above.

From the above discussion, every $\mathcal{A}$-path-type pseudo
tensor algebra $\mathcal{PT}(A,M)$ can be used to construct its
corresponding $\mathcal{A}$-pseudo path algebra
$PSE_k(\Delta,\mathcal{A})$; but, from this $\mathcal{A}$-pseudo
path algebra $PSE_k(\Delta,\mathcal{A})$, we can get uniquely the
free $\mathcal{A}$-path-type pseudo tensor algebra
$\mathcal{PT}(A,M^{F})$. Thus, we have:

\begin{lemma}
Every $\mathcal{A}$-path-type pseudo tensor algebra
$\mathcal{PT}(A,M)$ can be used to construct uniquely the free
$\mathcal{A}$-path-type pseudo tensor algebra
$\mathcal{PT}(A,M^{F})$. There is a surjective $k$-algebra
morphism $\pi$:
$\mathcal{PT}(A,M^{F})\rightarrow\mathcal{PT}(A,M)$ such that
$\pi(_{i}M^{F}_{j})=_{i}M_{j}$ for any $i,j\in I$.
\end{lemma}
{\em Proof}: It is sufficient to prove the second conclusion. For
$\mathcal{PT}(A,M)$, let the rank of $_{i}M_{j}$ be $r_{ij}$.
Thus, for the corresponding $\mathcal{A}$-pseudo path algebra
$PSE_k(\Delta,\mathcal{A})$, the number of the arrows from $i$ to
$j$ is $r_{ij}$, and then, in $\mathcal{PT}(A,M^{F})$,  the rank
of the free generators of $_{i}M^{F}_{j}$ given by the arrows is
also $r_{ij}$. Define $\pi$: $\mathcal{PT}(A,M^{F})\rightarrow
\mathcal{PT}(A,M)$ by giving a bijection between the set of the
free generators of $_{i}M^{F}_{j}$ and the set of the chosen
generators of $_{i}M_{j}$ with number of the rank. Then $\pi$ can
be expanded to become into a surjective $k$-algebra morphism with
$\pi(_{i}M^{F}_{j})=_{i}M_{j}$ for any $i,j\in I$.
$\;\;\;\;\square$

Next, we will show in the following Proposition 2.8 that every
$\mathcal{A}$-path-type pseudo tensor algebra is a homomorphic
image of its corresponding $\mathcal{A}$-pseudo path algebra.

The following criterion  for constructing algebra morphisms from
pseudo tensor algebras to other algebras is useful, which is
alternated from Lemma III.1.2 in \cite{ARS}. Contrast it with
Lemma 2.4.
\begin{lemma}
Let $A$ be a $k$-algebra and $M$ an $A$-bimodule. Let $\Lambda$ be
a $k$-algebra and $f: A\oplus M\rightarrow\Lambda$ a $k$-linear
map such that
 $f|_{A}: A\rightarrow\Lambda$ is an algebra morphism. Then there is a unique algebra homomorphism $\widetilde{f}:
 \mathcal{PT}(A,M)\rightarrow\Lambda$ such that  $\widetilde{f}|_{A\oplus
 M}=f$ and generally, $\widetilde{f}(\sum_{n=0}^{\infty}m^{n}_{1}\otimes_k \cdot\cdot\cdot\otimes_k m^{n}_{n})
 =\sum_{n=0}^{\infty}f(m^{n}_{1})\cdot\cdot\cdot f(m^{n}_{n})$ for $m^{n}_{1}\otimes_k\cdot\cdot\cdot\otimes_k m^{n}_{n}\in M(n)$.
\end{lemma}
{\em Proof:}$\;\;$ Consider the map $\phi: M\times M\rightarrow
\Lambda$ defined by $\phi(m_{1},m_{2})=f(m_{1})f(m_{2})$ for
$m_{1}$ and $m_{2}$ in $M$. We have for $\alpha\in k$ that
$\phi(m_{1}\alpha,m_{2})=f(m_{1}\alpha)f(m_{2})=f(m_{1})f(\alpha
m_{2})=\phi(m_{1},\alpha m_{2})$. Hence there is a unique group
morphism $f_{2}: M\otimes_{k}M\rightarrow\Lambda$ such that
$f_{2}(m_{1}\otimes_km_{2})=f(m_{1})f(m_{2})$.  Moreover, $f_{2}$
is a $k$-linear map. Similarly, for the map $\phi: M\times
A\rightarrow \Lambda$ defined by $\phi(m,a)=f(m)f(a)$ for $m\in M$
and $a\in A$, one can induce the $k$-linear map $f_{2}:
M\otimes_{k}A\rightarrow\Lambda$ satisfying $f_{2}(m\otimes_k
a)=f(m)f(a)$.

By induction, we can obtain the unique $k$-linear map $f_{n}:
M(n)\rightarrow\Lambda$ satisfying
$f_{n}(v_{1}\otimes_k\cdot\cdot\cdot\otimes_k
v_{n})=f(v_{1})\cdot\cdot\cdot f(v_{n})$. Since $f|_{A}$ is a
$k$-algebra homomorphism, we define $\widetilde{f}:
 \mathcal{PT}(A,M)\rightarrow\Lambda$ by  $\widetilde{f}|_{A\oplus
 M}=f$ and $\widetilde{f}(\sum_{n=0}^{\infty}m^{n}_{1}\otimes_k\cdot\cdot\cdot\otimes_km^{n}_{n})
 =\sum_{n=0}^{\infty}f(m^{n}_{1})\cdot\cdot\cdot f(m^{n}_{n})$ for $m^{n}_{1}\otimes_k\cdot\cdot\cdot\otimes_km^{n}_{n}\in M(n)$,
 which can be seen easily to be a $k$-algebra homomorphism uniquely determined by
 $f$.

 In fact, for $m_{1}\otimes_k\cdot\cdot\cdot\otimes_km_{n}\in
 M(n)$ and $\bar m_{1}\otimes_k\cdot\cdot\cdot\otimes_k\bar m_{l}\in
 M(l)$, if $m_{n}, \bar m_{1}\in A$, then $\widetilde{f}((m_{1}\otimes_k\cdot\cdot\cdot\otimes_km_{n})\cdot
 (\bar m_{1}\otimes_k\cdot\cdot\cdot\otimes_k\bar m_{l}))=\widetilde{f}(m_{1}\otimes_k\cdot\cdot\cdot\otimes_km_{n-1}\otimes_km_{n}\otimes_A
 \bar m_{1}\otimes_k\bar m_{2}\otimes_k\cdot\cdot\cdot\otimes_k\bar m_{l})\\
 =\widetilde{f}(m_{1}\otimes_k\cdot\cdot\cdot\otimes_km_{n-1}\otimes_km_{n}\bar m_{1}\otimes_k\bar m_{2}\otimes_k\cdot\cdot\cdot\otimes_k\bar m_{l})
 =f(m_{1})\cdot\cdot\cdot f(m_{n-1})f(m_{n}\bar m_{1})f(\bar m_{2})\cdot\cdot\cdot f(\bar
 m_{l})\\
 =f(m_{1})\cdot\cdot\cdot f(m_{n-1})f(m_{n})f(\bar m_{1})f(\bar m_{2})\cdot\cdot\cdot f(\bar
 m_{l})=\widetilde{f}(m_{1}\otimes_k\cdot\cdot\cdot\otimes_km_{n})\widetilde{f}
(\bar m_{1}\otimes_k\cdot\cdot\cdot\otimes_k\bar m_{l})$;\\
In the other cases, it can be proved similarly.$\>\>$
$\;\;\;\;\square$
\\

Comparing the definitions of generalized path algebra, tensor
algebra and pseudo path algebra, pseudo tensor algebra, the
following facts hold:
\begin{fact} {\em (1)}$\;$ There is a natural surjective
homomorphism$\;$ $\iota:
PSE_{k}(\Delta,\mathcal{A})\longrightarrow k(\Delta,\mathcal{A})$
with  $$ker\iota = \langle a\beta b\cdot c-a\beta bc,\; c\cdot
a\beta b-ca\beta b,\; a\alpha b\cdot c\beta d-a\alpha 1\cdot
bc\cdot 1\beta d\rangle$$ for any $a, b, c, d\in A=\oplus_i
A_{i}$, $\alpha, \beta\in\Delta_1$, where $1$ is the identity of
$A$. It follows that $$PSE_{k}(\Delta,\mathcal{A})/ker\iota\cong
k(\Delta,\mathcal{A})$$ as algebras.

 {\em (2)}$\;$ There is a
natural surjective homomorphism$\;$ $\tau:
\mathcal{PT}(A,M)\longrightarrow T(A,M)$ with $$ker\tau = \langle
m\otimes c-mc\otimes 1,\; c\otimes m-1\otimes cm,\; mb\otimes
cn-m\otimes bc\otimes n\rangle$$ for any $b, c\in A$, $m, n\in M$,
where $1$ is the identity of $A$. It follows that
$$\mathcal{PT}(A,M)/ker\tau\cong T(A,M)$$ as algebras.
\end{fact}

Clearly, all $\mathcal{A}$-pseudo-paths of length zero
(equivalently, $\mathcal{A}$-paths of length zero), that is, the
elements of $\bigcup_{i\in\Delta_{0}}A_{i}$ can generate a
subalgebra of $PSE_{k}(\Delta,\mathcal{A})$ (respectively,
$k(\Delta,\mathcal{A})$). Denote this subalgebra by
$PSE_{k}(\Delta_{0},\mathcal{A})$ (respectively,
$k(\Delta_{0},\mathcal{A})$). Then,
$PSE_{k}(\Delta_{0},\mathcal{A})=k(\Delta_{0},\mathcal{A})$, or
say, $\iota_{PSE_{k}(\Delta_{0},\mathcal{A})}=id$. Denote by
$PSE_{k}(\Delta_{1},\mathcal{A})$ (respectively,
$k(\Delta_{1},\mathcal{A})$) the $k$-linear space consisting of
all pure $\mathcal{A}$-pseudo-paths (respectively,  all
 $\mathcal{A}$-paths) of length $1$ and by $J$ (respectively,  $\widetilde{J}$) the ideal in
 $PSE_{k}(\Delta,\mathcal{A})$ (respectively,  $k(\Delta,\mathcal{A})$)
generated by all elements in $PSE_{k}(\Delta_{1},\mathcal{A})$
(respectively,  $k(\Delta_{1},\mathcal{A})$).

It is easy to see that $PSE_{k}(\Delta_{1},\mathcal{A})$
(respectively, $k(\Delta_{1},\mathcal{A})$) is an $A$-sub-bimodule
of $PSE_{k}(\Delta,\mathcal{A})$ (respectively,
$k(\Delta,\mathcal{A})$), and

 (i)
 $\;\;\iota(PSE_{k}(\Delta_1,\mathcal{A}))=k(\Delta_1,\mathcal{A})$;

 (ii) $\;\;\iota J=\widetilde{J}$, $\;\;\iota^{-1}\widetilde{J}=J$.

Now, we will show some useful properties of
$\mathcal{A}$-pseudo-path algebras which hold similarly for
$\mathcal{A}$-path
 algebra under the relationships in Fact 2.5.

\begin{lemma}
Let $\mathcal{PT}(A,M^{F})$ be the free $\mathcal{A}$-path-type
pseudo tensor algebra built by an $\mathcal{A}$-pseudo path
algebra $PSE_{k}(\Delta,\mathcal{A})$. Then there is a $k$-algebra
isomorphism $\phi$: $\mathcal{PT}(A,M^{F})\rightarrow
PSE_{k}(\Delta,\mathcal{A})$ such that for any $t\geq 1$,
$$\phi(\bigoplus_{n,l\geq
 t}M^{F}(n,l))=J^{t}.$$
\end{lemma}
{\em Proof}: By the multiplication in
$PSE_{k}(\Delta,\mathcal{A})$, $[a_{i}]\cdot[a_{j}]=0$ for
$i\not=j$ and $a_{i}\in A_{i}$, $a_{j}\in A_{j}$. Obviously, we
have a $k$-algebra isomorphism $f$: $A=\bigoplus_{i\in
I}A_{i}\rightarrow PSE_{k}(\Delta_{0},\mathcal{A})$ by
$f(a_{1}+\cdot\cdot\cdot+a_{n})=[a_{1}]+\cdot\cdot\cdot+[a_{n}]$.
And, define $f$: $M^{F}=\bigoplus_{i,j\in
I}$$_{i}M^{F}_{j}\rightarrow PSE_{k}(\Delta_{1},\mathcal{A})$ by
giving a bijection between a chosen basis for each $_{i}M^{F}_{j}$
and the set of arrows from $i$ to $j$, that is,
$f(am_{\alpha_{ij}}b)=a\alpha_{ij}b$ where $\alpha_{ij}$ is an
arrow from $i$ to $j$ and $m_{\alpha_{ij}}$ is the correspondent
element in the basis of $_{i}M^{F}_{j}$, $a,b\in A$.
 Since
 $PSE_{k}(\Delta_{0},\mathcal{A})$ is
 a $k$-subalgebra of
$PSE_{k}(\Delta,\mathcal{A})$, there is by Lemma 2.4 a $k$-algebra
morphism $\widetilde{f}$:
 $\mathcal{PT}(A,M^{F})\rightarrow PSE_{k}(\Delta,\mathcal{A})$ such that  $$\widetilde{f}|_{A\oplus
 M^{F}}=f$$
  and $$\widetilde{f}(\sum_{n=0}^{\infty}m^{n}_{1}\otimes\cdot\cdot\cdot\otimes m^{n}_{n})
 =\sum_{n=0}^{\infty}f(m^{n}_{1})\cdot...\cdot f(m^{n}_{n})$$\\
 for $m^{n}_{1}\otimes\cdot\cdot\cdot\otimes m^{n}_{n}\in M^{F}(n)$.
 Thus, $\widetilde{f}((A\otimes_{k}M^{F}\otimes_{k}A)^{t})=(A\cdot PSE_{k}(\Delta_{1},\mathcal{A})\cdot A)^{t}$
 and moreover, $\widetilde{f}(\bigoplus_{n,\;l\geq
 t}M^{F}(n,l))=J^{t}$, in particular, $\widetilde{f}(\bigoplus_{k\geq
 1}M^{F}(k))=J$. But,
 $PSE_{k}(\Delta,\mathcal{A})=PSE_{k}(\Delta_{0},\mathcal{A})\cup J\cup\cdot\cdot\cdot\cup
 J^{t}\cup\cdot\cdot\cdot$. Hence $\widetilde{f}$ is surjective.

 Let $\{x_{\lambda}\}$ denote a $k$-basis of $A$.
 For $M^{F}(n,l)$, we have a $k$-basis formed by some elements as
 $$x_{\lambda_{i_{1}}}\otimes
 x_{\lambda_{j_{1}}}m_{1}x_{\lambda_{k_{1}}}\otimes x_{\lambda_{i_{2}}}\otimes x_{\lambda_{j_{2}}}m_{2}x_{\lambda_{k_{2}}}
\otimes\cdot\cdot\cdot\otimes x_{\lambda_{i_{l}}}\otimes
x_{\lambda_{j_{l}}}m_{l}x_{\lambda_{k_{l}}}\otimes\cdot\cdot\cdot$$
 where there is some $\mathcal{A}$-pseudo-path
 $$[x_{\lambda_{i_{1}}}\cdot
 x_{\lambda_{j_{1}}}\beta_{1}x_{\lambda_{k_{1}}}\cdot x_{\lambda_{i_{2}}}\cdot x_{\lambda_{j_{2}}}\beta_{2}
 x_{\lambda_{k_{2}}}
\cdot...\cdot x_{\lambda_{i_{l}}}\cdot
x_{\lambda_{j_{l}}}\beta_{}x_{\lambda_{k_{l}}}\cdot...]$$
 in $PSE_{k}(\Delta,\mathcal{A})$ such that
$m_{j}$ is amongst the chosen basis elements in
$_{s(\beta_{j})}M^{F}_{s(\beta_{j+1})}$ of the correspondent arrow
$\beta_{j}$ for $j=1,...,t$.  Then

$\widetilde{f}(x_{\lambda_{i_{1}}}\otimes
 x_{\lambda_{j_{1}}}m_{1}x_{\lambda_{k_{1}}}\otimes x_{\lambda_{i_{2}}}\otimes x_{\lambda_{j_{2}}}m_{2}x_{\lambda_{k_{2}}}
\otimes\cdot\cdot\cdot\otimes x_{\lambda_{i_{t}}}\otimes
x_{\lambda_{j_{l}}}m_{l}x_{\lambda_{k_{l}}}\otimes\cdot\cdot\cdot)\\
=[x_{\lambda_{i_{1}}}\cdot
 x_{\lambda_{j_{1}}}\beta_{1}x_{\lambda_{k_{1}}}\cdot x_{\lambda_{i_{2}}}\cdot x_{\lambda_{j_{2}}}\beta_{2}
 x_{\lambda_{k_{2}}}
\cdot...\cdot x_{\lambda_{i_{l}}}\cdot
x_{\lambda_{j_{l}}}\beta_{l}x_{\lambda_{k_{l}}}\cdot...]$
\\
\\
It implies that distinct basis elements are mapped to distinct
$\mathcal{A}$-pseudo-paths.
 And, for $a_{1}+\cdot\cdot\cdot+a_{n}\not=0$ in $A$,
 $f(a_{1}+\cdot\cdot\cdot+a_{n})=[a_{1}]+\cdot\cdot\cdot+[a_{n}]\not=0$.
 Hence $\widetilde{f}$ is
injective. Therefore, $\phi=\widetilde{f}$ is a $k$-algebra
isomorphism with the desired properties.
$\;\;\;\;\square$
\\

By Lemma 2.6, $PSE_{k}(\Delta,\mathcal
A)\stackrel{\phi^{-1}}{\cong}\mathcal{PT}(A,M^F)$.  Then,
$ker\iota\stackrel{\phi^{-1}}{\cong}ker\tau$. Thus, a natural
induced algebra homomorphism  $\overline{\phi^{-1}}$ is obtained
from  $\phi^{-1}$ so that $PSE_{k}(\Delta,\mathcal
A)/ker\iota\stackrel{\overline{\phi^{-1}}}{\cong}\mathcal{PT}(A,M^F)/ker\tau$.
Moreover, by Fact 2.5, we get the following $\widetilde{\phi}$
from $\overline{\phi^{-1}}$ as above so as to gain the result on
$\mathcal{A}$-path algebra as similar as on
$\mathcal{A}$-pseudo-path algebra:
\begin{lemma}
Let $T(A,M^{F})$ be the free $\mathcal{A}$-path-type tensor
algebra built by a $\mathcal{A}$-path algebra
$k(\Delta,\mathcal{A})$. Then there is a $k$-algebra isomorphism
$\widetilde{\phi}$: $T(A,M^{F})\rightarrow k(\Delta,\mathcal{A})$
 such that for any $t\geq 1$,
$$\widetilde{\phi}(\bigoplus_{j\geq
 t}M^{Fj})=\widetilde J^{t}.$$
\end{lemma}

From this, we obtain the commutative diagram:
\begin{equation}
 \begin{diagram}
\node{\mathcal{PT}(A,M^F)}\arrow[2]{e,t}{\cong,\;\;\phi}\arrow{s,r}{\tau}\node[2]{PSE_k(\Delta,\mathcal A)}\arrow{s,r}{\iota}\\
\node{T(A,M^F)}\arrow[2]{e,t}{\cong,\;\;\widetilde{\phi}}\node[2]{k(\Delta,\mathcal
A)}
 \end{diagram}
 \end{equation}

\begin{proposition}
Let $\mathcal{PT}(A,M)$ be an $\mathcal{A}$-path-type pseudo
tensor algebra with the corresponding  $\mathcal{A}$-pseudo path
algebra $PSE_{k}(\Delta,\mathcal{A})$. Then there is a surjective
$k$-algebra homomorphism $\varphi$:
$PSE_{k}(\Delta,\mathcal{A})\rightarrow\mathcal{PT}(A,M)$ such
that for any $t\geq 1$,
$$\varphi(J^{t})=\bigoplus_{n,\;l\geq t}M(n,l).$$
\end{proposition}
{\em Proof}: Let $\mathcal{PT}(A,M^{F})$ be the free
$\mathcal{A}$-path-type pseudo tensor algebra built by the
$\mathcal{A}$-pseudo path algebra $PSE_{k}(\Delta,\mathcal{A})$.
Then by Lemma 2.6, there is a $k$-algebra isomorphism $\phi$:
$\mathcal{PT}(A,M^{F})\rightarrow PSE_{k}(\Delta,\mathcal{A})$
such that for any $t\geq 1$, $\phi(\bigoplus_{n,\;l\geq
 t}M^{F}(n,l))=J^{t}.$

 On the other hand, by Lemma 2.3,
there is a surjective $k$-algebra morphism $\pi$:
$\mathcal{PT}(A,M^{F})\rightarrow\mathcal{PT}(A,M)$ such that
$\pi(_{i}M^{F}_{j})=_{i}M_{j}$ for any $i,j\in I$, so
$\pi(M^{F})=M$.

  Therefore,
$\varphi=\pi\phi^{-1}$: $PSE_{k}(\Delta,\mathcal{A})\rightarrow
\mathcal{PT}(A,M)$ is a surjective $k$-algebra morphism with
$\varphi(J^{t})=\pi(\bigoplus_{n,\;l\geq
 t}M^{F}(n,l))=\bigoplus_{n,\;l\geq
 t}M(n,l)$
 for any $t\geq 1$.
$\;\;\;\;\square$
\\

According to $\varphi=\pi\phi^{-1}$ and the description of
ker$\iota$ and ker$\tau$ in Fact 2.5, we have
$\varphi($ker$\iota)=$ker$\tau$. Then, by Proposition 2.8, we
induce naturally a surjective $k$-algebra homomorphism

$\widetilde{\varphi}: PSE_k(\Delta,\mathcal
A)/ker\iota\rightarrow\varphi(PSE_k(\Delta,\mathcal A))/\varphi
(ker\iota)=\mathcal{PT}(A,M)/ker\tau.$\\
Thus, the similar result holds for $\mathcal A$-path-type tensor
algebra:
\begin{proposition}
Let $T(A,M)$ be an $\mathcal{A}$-path-type tensor algebra with the
corresponding  $\mathcal{A}$-path algebra $k(\Delta,\mathcal{A})$.
Then there is a surjective $k$-algebra homomorphism
$\widetilde{\varphi}$: $k(\Delta,\mathcal{A})\rightarrow T(A,M)$
such that for any $t\geq 1$,
$$\widetilde{\varphi}(\widetilde{J}^{t})=\bigoplus_{j\geq t}M^j.$$
\end{proposition}

Also, we obtain the commutative diagram:
\begin{equation}
 \begin{diagram}
\node{PSE_k(\Delta,\mathcal A)}\arrow[2]{e,t}{\varphi}\arrow{s,r}{\iota}\node[2]{\mathcal{PT}(A,M)}\arrow{s,r}{\tau}\\
\node{k(\Delta,\mathcal
A)}\arrow[2]{e,t}{\widetilde{\varphi}}\node[2]{T(A,M)}
 \end{diagram}
 \end{equation}

 A {\em
 relation} $\sigma$ on an  $\mathcal{A}$-pseudo path algebra
$PSE_{k}(\Delta,\mathcal{A})$ (respectively,  $\mathcal{A}$-path
algebra $k(\Delta,\mathcal{A})$) is a $k$-linear combination of
some general $\mathcal{A}$-pseudo paths (respectively,  some
$\mathcal{A}$-paths) $P_{i}$ with the same start vertex and the
same end vertex, that is,
$\sigma=k_{1}P_{1}+\cdot\cdot\cdot+k_{n}P_{n}$ with $k_{i}\in k$
and $s(P_{1})=\cdot\cdot\cdot=s(P_{n})$ and
$e(P_{1})=\cdot\cdot\cdot=e(P_{n})$.  If
$\rho=\{\sigma_{t}\}_{t\in T}$ is a set of relations on
$PSE_{k}(\Delta,\mathcal{A})$ (respectively,
$k(\Delta,\mathcal{A})$), the pair
$(PSE_{k}(\Delta,\mathcal{A}),\rho)$ (respectively,
$(k(\Delta,\mathcal{A}),\rho)$) is called an {\em
$\mathcal{A}$-pseudo path algebra with relations} (respectively,
{\em $\mathcal{A}$-path algebra with relations}). Associated with
$(PSE_{k}(\Delta,\mathcal{A}),\rho)$ (respectively,
$(k(\Delta,\mathcal{A}),\rho)$) is the quotient $k$-algebra
$PSE_{k}(\Delta,\mathcal{A},\rho)\stackrel{\rm
def}{=}PSE_{k}(\Delta,\mathcal{A})/\langle\rho\rangle$
(respectively, $k(\Delta,\mathcal{A},\rho)\stackrel{\rm
def}{=}k(\Delta,\mathcal{A})/\langle\rho\rangle$), where
$\langle\rho\rangle$ denotes the ideal in
$PSE_{k}(\Delta,\mathcal{A})$ (respectively,  in
$k(\Delta,\mathcal{A})$) generated by the set of relations $\rho$.
When the length $l(P_{i})$ of each $P_{i}$ is at least $j$, it
holds $\langle\rho\rangle\subset J^{j}$ (respectively,
$\langle\rho\rangle\subset \widetilde J^{j}$).

 For an element $x\in PSE_{k}(\Delta,\mathcal{A})$ (respectively,  $\in k(\Delta,\mathcal{A})$), we write by
$\bar{x}$ the corresponding element in
$PSE_{k}(\Delta,\mathcal{A},\rho)$ (respectively,
$k(\Delta,\mathcal{A},\rho)$).
\begin{fact}
$\delta\in k(\Delta,\mathcal A)$ is a relation if and only if all
$\sigma\in\iota^{-1}(\delta)$ are  relations on
$PSE_k(\Delta,\mathcal A)$.
\end{fact}

This fact can be seen easily from the definition of $\iota$. Note
that the lengths of paths in a relation do not be restricted here.
So, we have:
\begin{proposition}
Suppose that $\Delta$ is a finite quiver.

 {\em (i)}   Each element $x$ in
$PSE_{k}(\Delta,\mathcal{A})$ (respectively,
$k(\Delta,\mathcal{A})$) is a sum of some relations;

{\em (ii)}  Every ideal $I$ of $PSE_{k}(\Delta,\mathcal{A})$
(respectively, $k(\Delta,\mathcal{A})$) can be generated by a set
of relations.
\end{proposition}
{\em Proof}: (i) Let $1$ be the identity of $A$, $e_{i}$ the
identity of $A_{i}$ for $i\in\Delta_{0}$. Then
$1=\sum_{i\in\Delta}e_{i}$ is a decomposition into orthogonal
idempotents $e_{i}$.

$x=1\cdot x\cdot 1=\sum_{i,j\in\Delta_{0}}e_{i}\cdot x\cdot
e_{j}$. Due to the multiplication of
$\mid\Delta_{0}\mid=n<\infty$, $e_{i}\cdot x\cdot e_{j}$ can be
expanded as a $k$-linear combination of some such
$\mathcal{A}$-paths which have the same start vertex $i$ and the
same end vertex $j$, so
 $e_{i}\cdot x\cdot e_{j}$ is a relation on $PSE_{k}(\Delta,\mathcal{A})$.

 (ii) Assume $I$ is generated by
 $\{x_{\lambda}\}_{\lambda\in\Lambda}$. By (i), each $x_{\lambda}$
 is  a sum of some relations $\{\sigma_{\lambda,\;i}\}$. Then $I$
 is generated by all $\{\sigma_{\lambda,\;i}\}$.
 $\;\;\;\;\square$
\\

 By the definition of $J$, we have
$$PSE_{k}(\Delta,\mathcal{A},\rho)/\bar{J}
=(PSE_{k}(\Delta,\mathcal{A})/\langle\rho\rangle)/(J/\langle\rho\rangle)\cong
PSE_{k}(\Delta,\mathcal{A})/J\cong\oplus_{i\in\Delta_{0}}A_{i}.$$

Suppose all $A_{i}$ are $k$-simple algebras and
$J^{t}\subset\langle\rho\rangle$ for some integer $t$. Then
$PSE_{k}(\Delta,\mathcal{A},\rho)/\bar{J}\cong\oplus_{i\in\Delta_{0}}A_{i}$
is semisimple and $\bar{J}^{t}=0$. It follows that $\bar{J}={\rm
rad}PSE_{k}(\Delta,\mathcal{A},\rho)$. The similar discussion can
be taken for $\widetilde J$ of $k(\Delta,\mathcal A)$.  Hence we
get:
\begin{proposition}
{\em (i)} Let $(PSE_{k}(\Delta,\mathcal{A}),\rho)$ be an
$\mathcal{A}$-pseudo path algebra with relations where $A_{i}$ are
simple for all $i\in\Delta_{0}$. Assume that
$J^{t}\subset\langle\rho\rangle$ for some $t$. Then the image
$\bar{J}$ of $J$ in $PSE_{k}(\Delta,\mathcal{A},\rho)$ is ${\rm
rad}PSE_{k}(\Delta,\mathcal{A},\rho)$, that is,  $\bar{J}={\rm
rad}PSE_{k}(\Delta,\mathcal{A},\rho)$;

 {\em (ii)} Let
$(k(\Delta,\mathcal{A}),\rho)$ be an $\mathcal{A}$-path algebra
with relations where $A_{i}$ are simple for all $i\in\Delta_{0}$.
Assume that $\widetilde{J}^{t}\subset\langle\rho\rangle$ for some
$t$. Then the image $\overline{\widetilde{J}}$ of $\widetilde{J}$
in $k(\Delta,\mathcal{A},\rho)$ is ${\rm
rad}k(\Delta,\mathcal{A},\rho)$, that is,
$\overline{\widetilde{J}}={\rm rad}k(\Delta,\mathcal{A},\rho)$.

\end{proposition}

Now, suppose $A$ is a left Artinian algebra over $k$, $r=r(A)$ the
radical of $A$. Then for all $l\geq 0$, the ring $r^{l}/r^{l+1}$
is an $A$-bimodule by $a\cdot(r^{l}/r^{l+1})\cdot
b=ar^{l}b/r^{l+1}$ for $a,b\in A$. From $r\cdot r^{l}/r^{l+1}=0$
and $r^{l}/r^{l+1}\cdot r=0$, we know that $r^{l}/r^{l+1}$ is a
semisimple left and right $A$-module. For $\bar{x}=x+r\in A/r$,
let $\bar{x}\cdot(r^{l}/r^{l+1})\stackrel{\rm
def}{=}x\cdot(r^{l}/r^{l+1}) =xr^{l}/r^{l+1}$ and
$(r^{l}/r^{l+1})\cdot\bar{x}=(r^{l}/r^{l+1})\cdot
x=r^{l}x/r^{l+1}$, then $r^{l}/r^{l+1}$ is also an $A/r$-bimodule
and a semisimple left and right $A/r$-module.

\begin{proposition}\label{prop2.1}
Let $A$ be a left Artinian algebra over $k$, $r=r(A)$ the radical
of $A$. Write $A/r=\bigoplus_{i=1}^{s}\overline{A}_{i}$ where
$\overline{A}_{i}$ is a simple subalgebra for each $i$. Then,  for
all $l\geq 0$,

{\em (i)}$\;$ $r^{l}/r^{l+1}$ is finitely generated as an
$A/r$-bimodule;

{\em (ii)}$\;$ $_{i}M_{j}^{(l)}\stackrel{\rm
def}{=}\overline{A}_{i} \cdot r^{l}/r^{l+1}\cdot \overline{A}_{j}$
is finitely generated as
$\overline{A}_{i}$-$\overline{A}_{j}$-bimodule for each pair
$(i,j)$.
\end{proposition}
{\em Proof}: $\;$ (i)  Since $A$ is left Artinian, $r^{l}/r^{l+1}$
is finitely generated as a left $A$-module by Corollary I.3.2 in
\cite{ARS}. So, we can write
$r^{l}/r^{l+1}=\sum^{w}_{p=1}A\overline x_{p}$ with some
$\overline x_{p}\in r^{l}/r^{l+1}$.
 But, due to the definitions of actions, $A\overline x_{p}=(A/r)\overline x_{p}$
Then, $r^{l}/r^{l+1}=\sum^{w}_{p=1}(A/r)\overline x_{p}$.
Moreover, $r^{l}/r^{l+1}=r^{l}/r^{l+1}\cdot A/r
=(\sum^{w}_{p=1}(A/r)\overline
x_{p})(A/r)=\sum^{w}_{p=1}(A/r)\overline x_{p}(A/r)$,  which means
that $r^{l}/r^{l+1}$ is finitely generated as an $A/r$-bimodule.

(ii)  $_{i}M_{j}^{(l)}=\overline{A}_{i}\cdot r^{l}/r^{l+1}\cdot
\overline{A}_{j}=\overline{A}_{i}\cdot(\sum^{w}_{p=1}(A/r)\overline
x_{p}(A/r))\cdot\overline{A}_{j}=\sum^{w}_{p=1}\sum^{s}_{u,\;v=1}
\overline{A}_{i}\overline{A}_{u}\overline
x_{p}\overline{A}_{v}\overline{A}_{j}
=\sum^{w}_{p=1}\overline{A}_{i}\overline x_{p}\overline{A}_{j}$.
Hence, $_{i}M_{j}^{(l)}$ is finitely generated as
$\overline{A}_{i}$-$\overline{A}_{j}$-bimodule.
$\;\;\;\;\square$
\\

In particular, for $l=1$,  $\;_{i}M_{j}\stackrel{\rm
def}{=}\overline{A}_{i}\cdot r/r^{2}\cdot \overline{A}_{j}$ is
finitely generated as
$\overline{A}_{i}$-$\overline{A}_{j}$-bimodule for each pair
$(i,j)$. In the sequel, the rank of $_{i}M_{j}$ will be denoted by
$t_{ij}$.

For $k\not=i$, $\overline{A}_{k}\cdot
_{i}M_{j}=\overline{A}_{k}\cdot(\overline{A}_{i}\cdot r/r^{2}\cdot
\overline{A}_{j})=(\overline{A}_{k}\overline{A}_{i})\cdot(r/r^{2}\cdot
\overline{A}_{j})=0\cdot r/r^{2}\cdot\overline{A}_{j}=0$;
similarly, for $k\not=j$, $_{i}M_{j}\cdot\overline{A}_{k}=0$.
Thus, we obtain the $\mathcal{A}$-path-type pseudo tensor algebra
$\mathcal{PT}(A/r,r/r^{2})$, the $\mathcal{A}$-path-type tensor
algebra $T(A/r,r/r^{2})$ and the corresponding
$\mathcal{A}$-pseudo path algebra
 $PSE_{k}(\Delta,\mathcal{A})$ and $\mathcal{A}$-path algebra
 $k(\Delta,\mathcal{A})$,
with $\mathcal{A}=\{\overline{A}_{i}: i\in\Delta_{0}\}$, where
$\Delta$ is called as
 {\em the quiver of the left Artinian algebra $A$}.

In the following text,  $A$  is always a left Artinian algebra. We
will firstly show that under some important conditions, a left
Artinian algebra $A$ is isomorphic to some
$PSE_{k}(\Delta,\mathcal{A},\rho)$.

\section{When The Quotient Algebra Can Be Lifted}

Firstly, we introduce the concept of the set  of primitive
orthogonal simple subalgebras  of a left Artinian algebra. For a
left Artinian algebra $A$ and
$A/r=\bigoplus_{i=1}^{s}\overline{A}_{i}$ with  simple subalgebras
$\overline{A}_{i}$  for all $i$ where $r=r(A)$ is the radical of
$A$,
 assume there are simple  $k$-subalgebras  $B_{1}, \cdot\cdot\cdot, B_{s}$ of $A$ satisfying
$B_{i}\cong\overline{A}_{i}$ as $k$-algebras  for all $i$  under the canonical
 morphism $\eta$: $A\rightarrow A/r$ and
 $B_{i}B_{j}=\left\{\begin{array}{ll} B_{i},
  &  \mbox{if $i=j$}\\
0,  &  \mbox{if $i\not=j$}
\end{array}
\right. $.  Then, $\widehat{B}=\{B_{1}, \cdot\cdot\cdot, B_{s}\}$
is said to be {\em the set  of primitive orthogonal simple
subalgebras of $A$}.

Obviously, $\overline{A}_{i}\overline{A}_{j}
    =\left\{\begin{array}{ll}
\overline{A}_{i},
  &  \mbox{if $i=j$}\\
0,  &  \mbox{if $i\not=j$}
\end{array}
\right. $. By the definition,  $\eta(B_{i})=\overline{A}_{i}$ for
any $i$. Every $B_{i}$ is a simple $k$-subalgebra of $A$, so
$B=B_{1}+\cdot\cdot\cdot+B_{s}$ is a semisimple subalgebra of $A$.

Our original idea is to introduce the concept of primitive
orthogonal simple subalgebras as a generalization of primitive
orthogonal idempotents and then transplant the method of primitive
orthogonal idempotents in elementary algebras into a left Artinian
algebras.

In a left Artinian algebra $A$, we will show the existence of the
set of primitive orthogonal simple $k$-subalgebras when $A/r$ can
be lifted.

An algebra morphism $\varepsilon$: $A/r\rightarrow A$ satisfying
$\eta\varepsilon=1$ will be called a {\em lifting} of the quotient
algebra $A/r$. In this case, we call {\em $A/r$ can be lifted}.
Evidently, a lifting $\varepsilon$ is always a monomorphism and
Im$\varepsilon=B$ is a subalgebra of $A$ which is isomorphic to
$A/r$. Then, $B$ is semisimple.
 Moreover, $A=B\oplus r$ as a direct sum of $k$-linear
spaces. Hence, $A/r$ can be lifted if and only if $A$ is splitting
over its radical $r$.

Now, we assume that $A/r$ can be lifted such that $A=B\oplus r$ as
above. For the canonical morphism $\eta$: $A\rightarrow A/r$,
Im$\eta|_{B}=(B+r)/r=A/r$. And, Ker$\eta|_{B}=0$ due to $r\cap
B=0$.  Thus, $\eta(B)=A/r$ and $B\stackrel{\eta|_{B}}{\cong}A/r$
as $k$-algebras. Since $B$ is semisimple, we write
$B=\bigoplus_{i=1}^{s}B_{i}$ with simple $k$-subalgebras $B_{i}$
for all $i$. Then $B_{i}B_{j}=\left\{\begin{array}{ll} B_{i},
  &  \mbox{if $i=j$}\\
0,  &  \mbox{if $i\not=j$}
\end{array}
\right. $.
 Moreover,  $\eta(B)=\sum_{i=1}^{s}\eta(B_{i})$ where
$\eta(B_{i})$ are simple $k$-subalgebras of $A/r$. Denote
$\overline{A}_{i}=\eta(B_{i})$. Then, $\widehat{B}=\{B_{1},
\cdot\cdot\cdot, B_{s}\}$ is the set of primitive orthogonal
simple subalgebras of $A$.

\begin{lemma}
Assume $A$ is a left Artinian $k$-algebra with $r=r(A)$ the
radical of $A$ and $A/r$ can be lifted such that $A=B\oplus r$
with $\widehat{B}=\{B_{1},\cdot\cdot\cdot,B_{s}\}$ the set of
primitive orthogonal simple subalgebras of $A$ as constructed
above.
 Write $A/r=\bigoplus_{i=1}^{s}\overline{A}_{i}$ where
$\overline{A}_{i}$ are simple algebras for all $i$. The following
statements hold:

 {\em (i)} Let $\{r_{\lambda }: \lambda\in I\}$ be a set of elements in $r$
with the index set $I$ such that the images $\overline r_{\lambda
}$ in $r/r^{2}$ for all $\lambda\in I$ generate $r/r^{2}$ as an
$A/r$-bimodule. Then $B_{1}\cup\cdot\cdot\cdot\cup
B_{s}\cup\{r_{\lambda }: \lambda\in I\}$ generates $A$ as a
$k$-algebra;

{\em (ii)} There is a surjective $k$-algebra homomorphism
 $\widetilde{f}: \mathcal{PT}(A/r,r/r^{2})\rightarrow A$ with
$$\bigoplus_{n\geq
rl(A)}\;\bigoplus_{max\{rl(A),\; (n-1)/2\}\leq l\leq
n}M(n,l)\subset
 Ker\widetilde{f}\subset\bigoplus_{j\geq 2}M(j)$$\\
where rl$(A)$ denotes the Loewy length of $A$ as a left
$A$-module.
\end{lemma}
{\em Proof}: (i) Since $r$ is nilpotent, there is the least $m$
such that $r^{m}=0$ but $r^{m-1}\not=0$. It is easy to see that
$m$ is just the Loewy length rl$(A)$.

In the follows, we will prove this result by using induction on
$m$.

When $m=1$, then $r=0$ and $A$ is semisimple. Thus
$B_{i}=\overline{A}_{i}$. Hence $A$ is generated as a $k$-algebra
by $B_{1}\cup\cdot\cdot\cdot\cup B_{s}$.

When $m=2$, we have $r^{2}=0$. For the canonical morphism $\eta$,
$\eta(B_{i})=\overline{A}_{i}$. So, as a $k$-algebra, $A/r$ can be
generated by $(B_{1}+r)\cup\cdot\cdot\cdot\cup(B_{s}+r)$. Write
$A/r=\langle B_{1}+r,\cdot\cdot\cdot,B_{s}+r\rangle/r$.
 And, $\langle B_{1}+r,\cdot\cdot\cdot,B_{s}+r\rangle/r=(\langle
B_{1},\cdot\cdot\cdot,B_{s}\rangle+r)/r$. Thus, $A/r=(\langle
B_{1},\cdot\cdot\cdot,B_{s}\rangle+r)/r$. Hence $A=\langle
B_{1},\cdot\cdot\cdot,B_{s}\rangle+r$. But,
$r/r^{2}=\sum_{\lambda\in I}A/r\cdot\overline r_{\lambda}
=\sum_{\lambda\in I}A/r\cdot(r_{\lambda}+r^{2}) =\sum_{\lambda\in
I}(Ar_{\lambda}+r^{2})=(\sum_{\lambda\in I} Ar_{\lambda})+r^{2}$.
 Then from $r^{2}=0$, we get $r=\sum_{\lambda\in I}Ar_{\lambda}$.
 It follows $A=\langle
B_{1},\cdot\cdot\cdot,B_{s}\rangle+r=\langle
B_{1},\cdot\cdot\cdot,B_{s}\rangle+\sum_{\lambda\in I}(\langle
B_{1},\cdot\cdot\cdot,B_{s}\rangle+r)r_{\lambda}=\langle
B_{1},\cdot\cdot\cdot,B_{s}\rangle+\sum_{\lambda\in I}\langle
B_{1},\cdot\cdot\cdot,B_{s}\rangle r_{\lambda}=\langle
B_{1}\cup\cdot\cdot\cdot\cup B_{s}\cup\{r_{\lambda}: \lambda\in
I\}\rangle$ as a $k$-algebra.

Assume now that the claim holds for $m=l\geq 2$. Then consider the
claim in the case $m=l+1$.

Let $P$ be the $k$-subalgebra of $A$ generated by
$B_{1}\cup\cdot\cdot\cdot\cdot\cup B_{s}\cup\{r_{\lambda}:
\lambda\in I\}$. Firstly, we will show that $P/(P\cap
r^{l})=A/r^{l}$.

Since $(A/r^{l})/(r/r^{l})\cong A/r$ is semisimple,
$r(A/r^{l})=r/r^{l}$ holds. By the induction assumption,
$r^{l+1}=0$ and $r^{i}\not=0$ for any $i\leq l$. For any $t$,
$(r/r^{l})^{t}(A/r^{l})=r^{t}A/r^{l}=r^{t}/r^{l}$ since
$r^{t}A=r^{t}$ due to the existence of identity of $A$. Thus,
$(r/r^{l})^{t}(A/r^{l})=0$ if and only if $t\geq l$. (If there
were $t<l$ such that $r^{t}=r^{l}$, then $r^{t+1}=r^{l+1}=0$. It
contradicts to rl$(A)=m=l+1$). Therefore rl$(A/r^{l})=l$.

Let $\zeta$: $A\rightarrow A/r^{l}$ be the canonical morphism and
$\widetilde{B_{i}}=\zeta(B_{i})$ are simple algebras for all $i$,
$\pi$ the canonical morphism from $A/r^{l}$ to
$(A/r^{l})/(r/r^{l})=A/r$. Then $\pi\zeta=\eta$. It follows that
$\pi(\widetilde{B_{i}})=\overline{A}_{i}$. This means that
 $\widehat{\widetilde{B}}=\{\widetilde{B_{1}},\cdot\cdot\cdot,
\widetilde{B_{s}}\}$ is the set of primitive radical-orthogonal
simple algebras of $A/r^{l}$.
  We have that all elements in $\{\overline r_{\lambda}: \lambda\in I\}$ in
$r/r^{2}$ generate $r/r^{2}$ as an $A/r$-module. But,
$A/r\cong(A/r^{l})/(r/r^{l})$,
$r/r^{2}\cong(r/r^{l})/(r/r^{l})^{2}$. So, all elements in
$\{\overline r_{\lambda}: \lambda\in I\}$ in
$(r/r^{l})/(r/r^{l})^{2}$ generate $(r/r^{l})/(r/r^{l})^{2}$ as an
$(A/r^{l})/(r/r^{l})$-module. Let $\widetilde
r_{\lambda}=\zeta(r_{\lambda})\in r/r^{l}$. Then $\pi(\widetilde
r_{\lambda})=\overline r_{\lambda}$ .  Thus,  by the induction
assumption,
$\widetilde{B_{1}}\cup\cdot\cdot\cdot\cup\widetilde{B_{s}}
\cup\{\widetilde r_{\lambda}: \lambda\in I\}$ generates the
$k$-algebra $A/r^{l}$.

On the other hand, $B_{1}\cup\cdot\cdot\cdot\cdot\cup
B_{s}\cup\{r_{\lambda}: \lambda\in I\}$ generates $P$. Then
$\widetilde{B_{1}}\cup\cdot\cdot\cdot\cup\widetilde{B_{s}}
\cup\{\widetilde{r_{\lambda}}: \lambda\in I\}$ generates the
$k$-algebra $P/(P\cap r^{l})$. But, $P/(P\cap r^{l})$ can be
embedded into $A/r^{l}$. Therefore, we get that $P/(P\cap
r^{l})=A/r^{l}$.

Below it will be proved that in fact $P=A$, which means that
$B_{1}\cup\cdot\cdot\cdot\cup B_{s}\cup\{r_{\lambda}: \lambda\in
I\}$ generates $A$.

Let $x\in A$. Then there exists $y\in P$ such that
$x+r^{l}=y+P\cap r^{l}$. It follows $x-y\in r^{l}$. Thus there are
$\alpha_{i}\in r^{l-1}$ and $\beta_{i}\in r$ such that
$x-y=\sum^{q}_{i=1}\alpha_{i}\beta_{i}$. But, $\alpha_{i}+r^{l}$
and $\beta_{i}+r^{l}$ in $A/r^{l}$ and $A/r^{l}=P/(P\cap r^{l})$.
Then there are $a_{i}$ and $b_{i}$ in $P$ such that
$\alpha_{i}+r^{l}=a_{i}+P\cap r^{l}$ and
$\beta_{i}+r^{l}=b_{i}+P\cap r^{l}$. Due to $\alpha_{i}\in
r^{l-1}$ and $\beta_{i}\in r$, we have $a_{i}\in r^{l-1}$ and
$b_{i}\in r$. Let $a'_{i}=\alpha_{i}-a_{i}$ and
$b'_{i}=\beta_{i}-b_{i}$. Then $a'_{i},b'_{i}\in r^{l}$. Hence
$\alpha_{i}\beta_{i}=(a_{i}+a'_{i})(b_{i}+b'_{i})
=a_{i}b_{i}+a'_{i}b_{i}+a_{i}b'_{i}+a'_{i}b'_{i}=a_{i}b_{i}\in P$
for all $i$ where $a'_{i}b_{i}\in r^{l+1}=0$, $a_{i}b'_{i}\in
r^{2l-1}=0$, $a'_{i}b'_{i}\in r^{2l}=0$. It follows that $x-y\in
P$. Then $x\in P$.

(ii) $r/r^{2}=A/r\cdot r/r^{2}\cdot
A/r=\sum_{i,j=1}^{s}\overline{A}_{i}\cdot
r/r^{2}\cdot\overline{A}_{j}$ is a direct sum decomposition due to
$\overline{A}_{i}^{2}=\overline{A}_{i}$ and
$\overline{A}_{i}\overline{A}_{j}=0$ for $i\not=j$. Corresponding
to this, in $A$, we denote $W=\sum_{i,j=1}^{s}B_{i}rB_{j}$, where
$B_{i}\stackrel{\rm\eta}{\cong}\overline{A}_{i}$.  $W$ is a direct
sum of $B_{i}rB_{j}$ due to $B_{i}^{2}=B_{i}$ and $B_{i}B_{j}=0$
for $i\not=j$. Obviously, $W$ is a subalgebra of $r$ and then of
$A$. And, $r/r^{2}$ is a $(A/r)$-bimodule with the action of $A/r$
as above.

$(A/r)\oplus (r/r^{2})$ is a $k$-algebra in which the
multiplication is taken through that of $A/r$ and $r/r^{2}$ and
the $A/r$-bimodule action of $r/r^{2}$.

For each pair of integers $i,j$ with $1\leq i,j\leq s$, choose
elements $\{y_{u}^{ij}\}_{u\in \Omega_{ij} }$ in $B_{i}rB_{j}$
such that $\{\overline{y}_{u}^{ij}\}_{u\in \Omega_{ij} }$ is a
$k$-basis for $\overline{A}_{i}\cdot r/r^{2}\cdot\overline{A}_{j}$
for $\overline{y}_{u}^{ij}=y_{u}^{ij}+r^2$ the image in $r/r^{2}$.
Then $\bigcup_{i,j=1}^{s}\{\overline{y}_{u}^{ij}\}_{u\in
\Omega_{ij} }$ is a basis for $r/r^{2}$.
 It follows from (i)
that $\bigcup_{i,j,u}\{y_{u}^{i,j}\}_{u\in \Omega_{ij}}\cup
B_{1}\cup\cdot\cdot\cdot\cup B_{s}$ generates $A$ as a
$k$-algebra.

 It is easy to see that
$\{y_{u}^{ij}\}_{u\in \Omega_{ij} }$ is $k$-linear independent in
$B_{i}rB_{j}$. From the fact that $W$ is a direct sum of
$B_{i}rB_{j}$, it follows $\bigcup_{i,j=1}^{s}\{y_{u}^{ij}\}_{u\in
\Omega_{ij} }$ is a $k$-linear independent set in $W$.

Define $f: (A/r)\oplus (r/r^{2})\rightarrow A$ by
$f|_{\overline{A}_{i}}=\eta^{-1}$ and
$f(\overline{y}_{u}^{ij})=y_{u}^{ij}$. Then,  $f|_{A/r}$:
$A/r\rightarrow B=f(A/r)$ is a $k$-algebra isomorphism since
$B\stackrel{\eta|_{B}}{\cong}A/r$, and $f|_{r/r^{2}}$:
$r/r^{2}\rightarrow f(r/r^{2}) (\subset W\subset r)$ is an
isomorphism as $k$-linear spaces. Thus, $f: (A/r)\oplus
(r/r^{2})\rightarrow A$ is a $k$-linear map.
 Hence by Lemma 2.4, there is a unique
algebra morphism $\widetilde{f}$:
$\mathcal{PT}(A/r,r/r^{2})\rightarrow A$ such that
$\widetilde{f}|_{(A/r)\oplus(r/r^{2})}=f$. As said above,
$\bigcup_{i,j,u}\{y_{u}^{i,j}\}_{u\in \Omega_{ij}}\cup
B_{1}\cup\cdot\cdot\cdot\cup B_{s}$ generates $A$ as a
$k$-algebra. Therefore,  $\widetilde{f}$ is surjective.

By the definition of $\widetilde{f}$, we have
$\widetilde{f}((r/r^{2})^{j})=f(r/r^{2})^{j}\subset r^{j}\subset
r^{2}$ for $j\geq 2$, where $(r/r^{2})^{j}$ denotes
$r/r^{2}\otimes_{k}r/r^{2}\otimes_{k}\cdot\cdot\cdot\otimes_{k}r/r^{2}$
with $j$ copies of $r/r^{2}$. And, $f|_{A/r}$ and $f|_{r/r^{2}}$
are monomorphic. By the definition of $f$ on $A/r$ and $r/r^{2}$,
it is easy to see that $\widetilde{f}|_{(A/r)\oplus(r/r^{2})}$:
$(A/r)\oplus(r/r^{2})\rightarrow A$ is a monomorphism with image
intersecting $r^{2}$ trivially. As denoted in Section 2,
$M(n)=\sum_{M_{1},M_{2},\cdot\cdot\cdot,M_{n}}M_{1}\otimes_{k}M_{2}\otimes_{k}\cdot\cdot\cdot\otimes_{k}M_{n}$
 where $M_{i}$ is either $r/r^{2}$ (at least there exists one) or $A/r$
but no two $A/r$'s are neighbouring, then
$\mathcal{PT}(A/r,r/r^{2})=A/r\oplus M(1)\oplus
M(2)\oplus\cdot\cdot\cdot\oplus M(n)\oplus\cdot\cdot\cdot$.
 It follows that
Ker$\widetilde{f}\subset\bigoplus_{j\geq 2}M(j)$.

On the other hand, $M(n,l)$ equals the sum of these items
$M_{1}\otimes_{k}M_{2}\otimes_{k}\cdot\cdot\cdot\otimes_{k}M_{n}$
of $M(n)$ in which there are $l$ $M_{i}$'s equal to $r/r^{2}$  and
$M(n)=\sum_{(n-1)/2\leq l\leq n}M(n,l)$ as in Section 2.
 $\widetilde{f}((r/r^{2})^{j})=0$ for $j\geq
rl(A)$ since $r^{j}=0$ in this case, it follows
$\widetilde{f}(M(n,l))=0$ for any $n$ and possible $l\geq rl(A)$.
Therefore we get $$\bigoplus_{n\geq
rl(A)}\;\bigoplus_{max\{rl(A),\; (n-1)/2\}\leq l\leq
n}M(n,l)\subset
 Ker\widetilde{f}.\;\;\;\;\;\;\;\;\;\square$$

\begin{theorem}
{\em (Generalized Gabriel's Theorem Under Lifting)}

 Assume that
$A$ is a left Artinian $k$-algebra and $A/r$ can be lifted. Then,
 $A\cong PSE_{k}(\Delta,\mathcal{A},\rho)$ with
 $J^{s}\subset\langle\rho\rangle\subset
J$
 for some $s$, where $\Delta$ is the quiver of $A$
and $\rho$ is a set of relations on $PSE_{k}(\Delta,\mathcal{A})$.
\end{theorem}
 {\em Proof}:  Let $\Delta$ be the associated quiver of $A$. By Lemma 3.1(ii),
there is the surjective $k$-algebra morphism $\widetilde{f}$:
$\mathcal{PT}(A/r,r/r^{2})\rightarrow A$ with
$$\bigoplus_{n\geq
rl(A)}\;\bigoplus_{max\{rl(A),\; (n-1)/2\}\leq l\leq
n}M(n,l)\subset
 Ker\widetilde{f}\subset\bigoplus_{j\geq 2}M(j).$$\\
By Proposition 2.8, there is the surjective $k$-algebra
homomorphism $\varphi$:
$PSE_{k}(\Delta,\mathcal{A})\rightarrow\mathcal{PT}(A/r,r/r^{2})$
such that for any $t\geq 1$, $\varphi(J^{t})=\bigoplus_{n,l\geq
 t}M(n,l)$. Then
$\widetilde{f}\varphi$: $PSE_k(\Delta,\mathcal{A})\rightarrow A$
is a surjective $k$-algebra morphism with the kernel
$I=$Ker$(\widetilde{f}\varphi)=\varphi^{-1}($Ker$\widetilde{f})$.

But, $\varphi(J^{rl(A)})=\bigoplus_{n,\;l\geq
 rl(A)}M(n,l)$ and $\varphi(J^{2})=\bigoplus_{n,\;l\geq
 2}M(n,l)$. So, by Lemma 3.1(ii),
$\varphi(J^{rl(A)})\subset$Ker$\widetilde{f}\subset\varphi(J^{2})+M(2,1)+M(3,1)$.

One can show $J^{t}\subset\varphi^{-1}\varphi(J^{t})\subset
J^{t}+\phi(\bigoplus_{n}\bigoplus_{l\leq
 t-1}M^{F}(n,l))\cap\phi($Ker$\pi)$ for $t\geq
1$. In fact, trivially, $J^{t}\subset\varphi^{-1}\varphi(J^{t})$.
On the other hand, $\varphi=\pi\phi^{-1}$ and then
$\varphi^{-1}=\phi\pi^{-1}$. By Proposition 2.8,
$\varphi(J^{t})=\bigoplus_{n,\;l\geq
 t}M(n,l)$. From the definition of $\pi$ in Lemma 2.3, it
can be seen that $\pi^{-1}(\bigoplus_{n,\;l\geq
 t}M(n,l))\subset\bigoplus_{n,\;l\geq
 t}M^{F}(n,l)+(\bigoplus_{n}\bigoplus_{l\leq
 t-1}M^{F}(n,l))\cap$Ker$\pi$. Thus, by
Lemma 2.3, we have

 $\varphi^{-1}\varphi(J^{t})=\phi\pi^{-1}(\bigoplus_{n,\;l\geq
 t}M(n,l))\subset\phi(\bigoplus_{n,\;l\geq
 t}M^{F}(n,l))+\phi(\bigoplus_{n}\bigoplus_{l\leq
 t-1}M^{F}(n,l))\cap\phi($Ker$\pi)\\
=J^{t}+\phi(\bigoplus_{n}\bigoplus_{l\leq
 t-1}M^{F}(n,l))\cap\phi($Ker$\pi)$. \\
  Hence,

$J^{rl(A)}\subset\varphi^{-1}\varphi(J^{rl(A)})
\subset\varphi^{-1}($Ker$\widetilde{f})=I\subset\varphi^{-1}\varphi(J^{2})+\varphi^{-1}(M(2,1)+M(3,1))\\
\subset
J^{2}+\phi(M^{F}(3,1)+M^{F}(2,1)+M^{F}(1,1))\cap\phi($Ker$\pi)+\varphi^{-1}(M(2,1)+M(3,1))
\\=J^{2}+A\cdot PSE(\Delta_{1},\mathcal{A})\cdot A$\\
since $\phi(M^{F}(1,1))\cap\phi($Ker$\pi)=0$ and then
$\phi(M^{F}(3,1)+M^{F}(2,1)+M^{F}(1,1))\cap\phi($Ker$\pi)+\varphi^{-1}(M(2,1)+M(3,1))
=A\cdot PSE(\Delta_{1},\mathcal{A})\cdot A$. \\
But, it is clear that $J^{2}+A\cdot
PSE(\Delta_{1},\mathcal{A})\cdot A = J$. Therefore, we get:
$$J^{rl(A)}\subset\varphi^{-1}(Ker\widetilde{f})=I\subset J$$

Lastly, by Proposition 2.11, there is a set $\rho$ of relations so
that $I$ can be generated by $\rho$, that is,
$I=\langle\rho\rangle$. Hence,
$PSE_k(\Delta,\mathcal{A},\rho)=PSE_k(\Delta,\mathcal{A})/\langle\rho\rangle\cong
A$
  with
  $\langle\rho\rangle=$Ker$(\widetilde{f}\varphi)$ and $J^{rl(A)}\subset\langle\rho\rangle
  \subset J$.
$\;\;\;\;\square$
\\

  Usually, for a left Artinian algebra $A$, the set $\rho$ of relations in Theorem
  3.2 is infinite. But, when $A$ is finite dimensional, we can
  show $\rho$ is finite.

  In fact, suppose that $A$ is finite dimensional, then $\overline
  A_{i}$ is finite dimensional for all $i$. Thus, the $k$-space
  consisting of all $\mathcal A$-pseudo paths of a certain length
  is finite dimensional. It follows that $J^{rl(A)}$ is the ideal
  finitely generated in $PSE_{k}(\Delta,\mathcal A)$ by all $\mathcal A$-pseudo paths of  length
 $rl(A)$. Similarly, $PSE_{k}(\Delta,\mathcal A)/J^{rl(A)}$ is
 generated finitely as a $k$-space, under the meaning of
 isomorphism, by all $\mathcal A$-paths of length less than
 $rl(A)$, so as well as $I/J^{rl(A)}$ as a $k$-subspace. Then it
 is easy to know that $I$ is a finitely generated ideal in
 $PSE_{k}(\Delta,\mathcal A)$. Assume
$\{\sigma_{1},\cdot\cdot\cdot,\sigma_{p}\}$ is a set of finite
generators for the ideal $I$. For the identity $\overline{1}$ of
$A/r$, we have the decomposition of orthogonal idempotents
$\overline{1}=\overline{e}_{1}+\cdot\cdot\cdot+\overline{e}_{s}$,
where $\overline{e}_{i}$ is the identity of $\overline{A}_{i}$.
Then
$\sigma_{l}=\overline{1}\cdot\sigma_{l}\cdot\overline{1}=\sum_{1\leq
i,j\leq s}\overline{e}_{i}\cdot\sigma_{l}\cdot\overline{e}_{j}$,
where  $\overline{e}_{i}\sigma_{l}\overline{e}_{j}$ can be
expanded as a
 $k$-linear combination of some such $\mathcal{A}$-pseudo paths which
 have the same start vertex $i$ and the same end vertex $j$. So,
 $\sigma^{ilj}=\overline{e}_{i}\cdot\sigma_{l}\cdot\overline{e}_{j}$
 is a relation on the $\mathcal{A}$-pseudo path algebra
 $PSE_{k}(\Delta,\mathcal{A})$. Moreover, $I$ is generated by all $\sigma^{ilj}$
 due to $\sigma_{l}=\sum_{i,j}\sigma^{ilj}$. Therefore, we have a
 finite set $\rho=\{\sigma^{ilj}: 1\leq i,j\leq s, 1\leq l\leq p\}$ with $I=\langle\rho\rangle$ such that $PSE_k(\Delta,\mathcal{A},\rho)=PSE_k(\Delta,\mathcal{A})/\langle\rho\rangle\cong
  A$. Therefore, the following holds:

  \begin{corollary}
 Assume that
$A$ is a finite dimensional $k$-algebra and $A/r$ can be lifted.
Then,
 $A\cong PSE_{k}(\Delta,\mathcal{A},\rho)$ with
 $J^{s}\subset\langle\rho\rangle\subset
J$
 for some $s$, where $\Delta$ is the quiver of $A$
and $\rho$ is a finite set of relations on
$PSE_{k}(\Delta,\mathcal{A})$.
  \end{corollary}

When $A$ is elementary, $A_{i}=A_{j}=k$ and $_{i}M_{j}=r/r^{2}$ is
free as a $k$-linear space. Thus, $\pi$ is an isomorphism, then
Ker$\pi=0$ and Ker$\varphi=0$. According to the classical Gabriel
  Theorem, we have $J^{rl(A)}\subset\langle\rho\rangle
  \subset J^{2}$, which is a special case of the results of Theorem 3.2 and Corollary 3.3.

  By the famous Wedderburn-Malcev Theorem (see \cite{P}), for a left Artinian $k$-algebra $A$ and its radical $r$,
    if Dim$A/r\leq 1$, then $A/r$ can be lifted. Here, the
    dimension Dim$A$ of a $k$-algebra $A$ is defined as
    \textrm{Dim}$A=\textrm{sup}\{n: H^{n}_{k}(A,M)\not=0$ for some $A$-bimodule
    $M\}$ where $H^{n}_{k}(A,M)$ means the $n$'th Hochschild
    cohomology module of $A$ with coefficients in $M$. In
    particular, Dim$A/r=0$ if and only if $A/r$ is a separable
    $k$-algebra. By Corollary 10.7b in \cite{P}, when $k$ is a
    perfect field (for example,   char$k=0$ or $k$ is a finite field), $A$
    is separable. So, we have the following:

\begin{proposition}
 Assume that
$A$ is a left Artinian $k$-algebra. Then,
 $A\cong PSE_{k}(\Delta,\mathcal{A},\rho)$ with $J^{s}\subset\langle\rho\rangle\subset J$
  for some $s$, where
$\Delta$ is the quiver of $A$ and $\rho$ is a set of relations of
$PSE_{k}(\Delta,\mathcal{A})$, if one of the following conditions
holds:

 {\em (i)}$\;\;$ {\em Dim}$A/r\leq 1$ for the radical $r$ of $A$;

 {\em (ii)}$\;\;$ $A/r$ is separable;

{\em (iii)}$\;\;$ $k$ is a perfect field (for example,   when {\em
char}$k=0$ or $k$ is a finite field).
\end{proposition}

Note that in Proposition 3.4,  the condition (ii) is a special
case of (i), and (iii) is that of (ii).

In Theorem 3.2, $A\cong PSE_{k}(\Delta,\mathcal{A},\rho)$ holds
where $\Delta$ is the quiver of $A$ from the corresponding
$\mathcal{A}$-pseudo path algebra of the $\mathcal{A}$-path-type
 pseudo tensor algebra $\mathcal{PT}(A/r,r/r^{2})$ by the
definitions in Section 2.
 Moreover, in the case that $\langle\rho\rangle\subset J_{\Delta}^2$, we will discuss the uniqueness of the correspondent pseudo path algebra and quiver of a left Artinian algebra
  up to isomorphism, that is, if there exists another quiver and its
 related pseudo path algebra so that the same isomorphism relation is satisfied.
 In fact, we have the
following statement on the uniqueness:
\begin{theorem}
Assume $A$ is a left Artinian  $k$-algebra. Let
$A/r(A)=\bigoplus^{p}_{i=1}\overline{A}_{i}$ with simple algebras
$\overline{A}_{i}$. If there is a quiver $\Delta$ and a pseudo
path algebra $PSE_{k}(\Delta,\mathcal{B})$ with a set of simple
algebras $\mathcal{B}=\{B_{1},\cdot\cdot\cdot,B_{q}\}$ and $\rho$
a set of relations satisfying that $A\cong
PSE_{k}(\Delta,\mathcal{B},\rho)$ with
$J_{\Delta}^{t}\subset\langle\rho\rangle\subset J_{\Delta}^2$ for
some $t$ and $J_{\Delta}$ the ideal in
$PSE_{k}(\Delta,\mathcal{B})$ generated by all pure paths in
$PSE_{k}(\Delta_{1},\mathcal{B})$, then $\Delta$ is just the
quiver of $A$ and $p=q$ such that $\overline{A}_{i}\cong B_{i}$
for $i=1,...,p$ after reindexing.
\end{theorem}
{\em Proof}:
$PSE_{k}(\Delta,\mathcal{B})/J_{\Delta}=B_{1}+\cdot\cdot\cdot+B_{q}$
due to the definition of $J_{\Delta}$. Since
$J_{\Delta}^{t}\subset\langle\rho\rangle$, it follows that
$(J_{\Delta}/\langle\rho\rangle)^{t}
=J_{\Delta}^{t}/\langle\rho\rangle=0$. And,

$PSE_{k}(\Delta,\mathcal{B},\rho)/(J_{\Delta}/\langle\rho\rangle)
 =(PSE_{k}(\Delta,\mathcal{B})/\langle\rho\rangle)/(J_{\Delta}/\langle\rho\rangle)
 =PSE_{k}(\Delta,\mathcal{B})/J_{\Delta}=B_{1}+\cdot\cdot\cdot+B_{q}$\\
  is semisimple. Hence,
$J_{\Delta}/\langle\rho\rangle$ is the radical of
$PSE_{k}(\Delta,\mathcal{B},\rho)$. Thus, from  $A\cong PSE
k(\Delta,\mathcal{B},\rho)$, it follows that $A/r(A)\cong
PSE_{k}(\Delta,\mathcal{B})/J_{\Delta}$. But,
 $A/r(A)=\bigoplus^{p}_{i=1}\overline{A}_{i}$ and
 $PSE_{k}(\Delta,\mathcal{B})/J_{\Delta}=B_{1}+\cdot\cdot\cdot+B_{q}$with
 $\overline{A}_{i}$ and $B_{j}$ are simple algebras. Therefore,
$p=q$ and $\overline{A}_{i}\cong B_{i}$ for $i=1,...,p$ after
reindexing, according to Wedderburn-Artinian Theorem.

On the other hand,  $A/r(A)^{2}\cong
PSE_{k}(\Delta,\mathcal{B})/J_{\Delta}^{2}$. Then, the quivers of
$A/r(A)^{2}$ and $PSE_{k}(\Delta,\mathcal{B})/J_{\Delta}^{2}$ are
the same.

But,
$PSE_{k}(\Delta,\mathcal{B})/J_{\Delta}^{2}=(PSE_{k}(\Delta,\mathcal{B})/\langle\rho\rangle)/(J_{\Delta}^{2}/\langle\rho\rangle)
=PSE_{k}(\Delta,\mathcal{B},\rho)/(J_{\Delta}^{2}/\langle\rho\rangle)$
and the radical of $PSE_{k}(\Delta,\mathcal{B},\rho)$ is
$J_{\Delta}/\langle\rho\rangle$. Then the radical of
$PSE_{k}(\Delta,\mathcal{B})/J_{\Delta}^{2}$ is
$(J_{\Delta}/\langle\rho\rangle)/(J_{\Delta}^{2}/\langle\rho\rangle)\cong
J_{\Delta}/J_{\Delta}^{2}$. So, the quivers of
 $PSE_{k}(\Delta,\mathcal{B})/J_{\Delta}^{2}$ are that of the $\mathcal{A}$-path-type pseudo tensor algebra
 $$\mathcal{PT}((PSE_{k}(\Delta,\mathcal{B})/J_{\Delta}^{2})/(J_{\Delta}/J_{\Delta}^{2}),\;
 (J_{\Delta}/J_{\Delta}^{2})/(J_{\Delta}^{2}/J_{\Delta}^{2}))\cong
 \mathcal{PT}(PSE_{k}(\Delta,\mathcal{B})/J_{\Delta},
 J_{\Delta}/J_{\Delta}^{2}).$$
Now, we consider the quiver $\Gamma$ of
$\mathcal{PT}(PSE_{k}(\Delta,\mathcal{B})/J_{\Delta},
 J_{\Delta}/J_{\Delta}^{2})$. From the definition of the quiver
 associating with an $\mathcal{A}$-path-type pseudo tensor algebra in Section
 2, we know that $\Gamma_{0}=\{1,\cdot\cdot\cdot,q\}=\Delta_{0}$.
 For any $i,j\in\Gamma_{0}$,  the number of arrows from $i$ to $j$ in
 $\Gamma$ is the rank $r_{ij}$ of $_{i}M_{j}=B_{i}\cdot J_{\Delta}/J_{\Delta}^{2}\cdot
 B_{j}$
 as a finitely generated $B_{i}$-$B_{j}$-bimodule. However, by the
 definition of $J_{\Delta}$, under the meaning of isomorphism, $B_{i}\cdot J_{\Delta}/J_{\Delta}^{2}\cdot
 B_{j}$ can be constructed as an $B_{i}$-$B_{j}$-linear expansion
 of all $\mathcal{A}$-pseudo-paths of length $1$ from $i$ to $j$ in
 $PSE_{k}(\Delta_{1},\mathcal{B})$. It means that $r_{ij}$ is equal to the number of arrows from $i$ to $j$ in
 $\Delta$. Thus, the number of arrows from $i$ to $j$ in
 $\Gamma$ is equal to that of arrows from $i$ to $j$ in
 $\Delta$. Then $\Gamma_{1}=\Delta_{1}$. Therefore, we get
 $\Gamma=\Delta$.

 Due to the above discussion, it implies that the quivers of
$A/r(A)^{2}$ is just $\Delta$.
 Moreover, $A/r(A)=(A/r(A)^{2})/(r(A)/r(A)^{2})$ and
 $r(A)/r(A)^{2}=(r(A)/r(A)^{2})/(r(A)/r(A)^{2})^{2}$, where $r(A)/r(A)^{2}$ is the radical of $A/r(A)^{2}$.
 So, the quivers $\Delta$ of
 $A/r(A)^{2}$ is also that of
 $$\mathcal{PT}((A/r(A)^{2})/(r(A)/r(A)^{2}),\; (r(A)/r(A)^{2})/(r(A)/r(A)^{2})^{2}).$$
But,
 $$\mathcal{PT}(A/r(A), r(A)/r(A)^{2})\cong\mathcal{PT}((A/r(A)^{2})/(r(A)/r(A)^{2}),\; (r(A)/r(A)^{2})/(r(A)/r(A)^{2})^{2}).$$
 It follows that $\Delta$ is the quiver of $A$.
$\;\;\;\;\square$
\\

 According to this theorem, we see that for a left Artinian algebra $A$, the existence of the pseudo path
 algebra such that $A$ is isomorphic to its quotient algebra (that is, Theorem 3.2) can deduce its
 uniqueness. That is, it is just the pseudo path algebra decided by the quiver
  and the semisimple decomposition of $A$.

  In Theorem 3.5, the condition $J_{\Delta}^{t}\subset\langle\rho\rangle\subset
  J_{\Delta}^2$ is usually said that the ideal $\langle\rho\rangle$ is
  {\em admissible}. Notice that in the conclusions of Theorem 3.2, Corollary 3.3 and Proposition 3.4, we
  have only $J_{\Delta}^{t}\subset\langle\rho\rangle\subset
\dot{}  J_{\Delta}$. So far, we can only say in these above
results, $\langle\rho\rangle$ seems not
  admissible, as well as the following Theorem 4.2, Corollary 4.3 and Proposition
  4.5.
  Therefore, we have not known whether the correspondent pseudo path algebra (respectively, the generalized
  path algebra) and the quiver of the algebra are unique up to
  isomorphism in these main results of Section 3 and 4.

Our main result, Theorem 3.2,  means when the quotient algebra of
a left Artinian algebra is lifted, the algebra can be covered by a
pseudo path algebra under an algebra homomorphism. We know that a
generalized path algebra must be a homomorphic image of a pseudo
path algebra and its definition seems to be more concise than that
of pseudo path algebra. So, it is natural to ask why we do not
look for a generalized path algebra to cover a left Artinian
algebra. In fact, this is our original idea. However,
unfortunately, in general, as shown by the following
counter-example, a left Artinian algebra with lifted quotient may
not be a
 homomorphic image of its correspondent $\mathcal A$-path-type tensor algebra. Thus, one cannot use the
 method as above (that is,  through Proposition 2.9) to gain  a generalized path algebra in order to cover the left Artinian algebra.
  The following counter-example was given by W.
Crawley-Boevey.
\begin{example}
There is an example of a finite dimensional algebra $A$ over a
field $k$ such that\\
{\em (a)} $A$ is splitting over its radical $r$, that is, $A/r$
can be
lifted;\\
 {\em (b)} there is {\bf no} surjective algebra homomorphism from
 $T(A/r,r/r^2)$ to $A$, that is, $A$ cannot be equivalent to some
 quotient of $T(A/r,r/r^2)$.
\end{example}

Concretely, we describe $A$ as follows:
\\
 (1)$\;$ Let $F/k$ be a finite field extension, and
let $\delta: F\rightarrow F$ be a non-zero $k$-derivation. For
example, one can take $k=\mathbf{Z}_2(t)$,
$F=\mathbf{Z}_2(\sqrt{t})$ and $\delta(p+q\sqrt{t})=q$ for $p,q\in
\mathbf{Z}_2(t)$ where $\mathbf{Z}_2$ denotes the prime field of
characteristic $2$. It is easy to check $\delta$ as a
$k$-derivation due to  $char k=2$.
\\
\\
(2)$\;$ Define $E=F\oplus F$ and consider it as an
$F$-$F$-bimodule with the actions:
$$f(x,y)=(fx,fy)\;\;\;\;\;\;\;\;(x,y)f=(xf+y\delta(f),yf)$$
for $x,y,f\in F$. Let $\theta$ and $\phi$ be $F$-$F$-bimodule
homomorphisms respectively from $F$ to $E$ and from $E$ to $F$
 satisfying
$$\theta(x)=(x,0)\;\;\;\;\;\;\;\phi(x,y)=y$$
for $x,y\in F$. Then we have the non-splitting extension of
$F$-$F$-bimodules:
$$0\rightarrow F\stackrel{\theta}{\rightarrow}E\stackrel{\phi}{\rightarrow}F\rightarrow 0$$

In fact, if there were $\psi: E\rightarrow F$ an $F$-$F$-bimodule
homomorphism with $\psi\cdot\theta=1_F$, then for all $f\in F$,
$$\delta(f)=\psi\theta(\delta(f))=\psi(\delta(f),0)=\psi(\delta,f)-\psi(0,f)
=\psi((0,1)f)-\psi(f(0,1))=\psi(0,1)f-f\psi(0,1)=0,$$ it follows
that $\delta=0$, which contradicts to the presumption on $\delta$.
\\
\\
 (3)$\;$  Define $A=F\oplus F\oplus E$ with multiplication given
by $$(x,y,e)(x',y',e')=(xx',xy'+yx',xe'+\theta(yy')+ex').$$

Let $S=\{(x,0,0): x\in F\}$. Then $S$ is a subalgebra of $A$
isomorphic to $F$.

Let $r=\{(0,y,e): y\in F, e\in E\}$. Then $r$ is an ideal in $A$
with $r^2=\{(0,0,e): e\in Im(\theta)\}$ and $r^3=0$. Thus $r$ is
the radical of $A$, and $A=S\oplus r$, so $A$ is splitting over
$r$.
\\
\\
(4)$\;$ As an $F$-$F$-module, $r/r^2$ is isomorphic to $F\oplus F$
due to the surjective $F$-$F$-module homomorphism $\pi:
r\rightarrow F\oplus F$ satisfying $\pi(0,y,e)=(y,\phi(e))$ with
ker$\pi=r^2$.
\\
\\
(5)$\;$ By (3) and (4), the $\mathcal A$-path-type tensor algebra
$T(A/r,r/r^2)\cong T(F,F\oplus F)$. Let $s=(1,0)$ and $t=(0,1)$,
then $F\oplus F\cong Fs\oplus Ft$. Thus, $T(F,F\oplus F)$
(equivalently, say $T(A/r,r/r^2)$) can be considered as the free
associative algebra $F\langle s,t\rangle$ generated by two
variables $s, t$ over $F$. It follows that the center
$Z(T(A/r,r/r^2))$ of $T(A/r,r/r^2)$ is equal to $F$.
 \\
\\
(6)$\;$  If $(x,y,e)\in Z(A)$ the center of $A$, then for all
$e'\in E$, $(x,y,e)$ commutes with $(0,0,e')$, thus
$(0,0,xe')=(0,0,e'x)$, so $xe'=e'x$. Taking $e'=(0,1)$, we get
$x(0,1)=(0,1)x$. But by (2), $x(0,1)=(0,x)$ and
$(0,1)x=(\delta(x),x)$. It follows that $\delta(x)=0$. Therefore,
we have:
$$Z(A)\subset\{(x,y,e): x,y\in F, e\in E, \delta(x)=0\}.$$
 (7)$\;$ If $L$ is a subalgebra of $Z(A)$ and is a field, then
 dim$_kL\stackrel{<}{\not=}$dim$_kF$.

 In fact, the composition $L\hookrightarrow Z(A)\hookrightarrow\{(x,y,e): x,y\in F, e\in E,
 \delta(x)=0\}\rightarrow\{x: \delta(x)=0\}$ is an algebra
 homomorphism. Assume $l=(x,y,e)\in L$ is in the kernel of this
 composition, then $x=0$ and $l=(0,y,e)$, so $l\in r$ the radical
 of $A$. By (3), $l^3=0$. But $L$ is a field, so $l=0$. It means
 that this composition is one-one map. Therefore,
$$dim_kL\leq dim_k\{x: \delta(x)=0\}\stackrel{<}{\not=}dim_kF$$
where ``$\not=$" is from $\delta\not=0$.
\\
\\
(8)$\;$ If there were a surjective algebra homomorphism $\lambda:
T(A/r,r/r^2)\rightarrow A$, it would induce a homomorphism of the
center $Z(T(A/r,r/r^2))$ of $T(A/r,r/r^2)$ into the center $Z(A)$
of $A$. By (5), $Z(T(A/r,r/r^2))=F$. Thus, $L=\lambda(F)$ would be
 a field and a subalgebra of $Z(A)$. By (7), we have
 $dim_kL\stackrel{<}{\not=}dim_kF$. On the other hand,
 if there is $x$ satisfying $0\not=x\in$Ker$\lambda|_F$, that is, $\lambda(x)=0$. Since $F$ is a
 field, we get $\lambda(1)=\lambda(1/x)\lambda(x)=0$, then
 $\lambda=0$ is induced due to $\lambda$ as an algebra
 homomorphism. It is impossible since $\lambda$ is surjective. It
 means Ker$\lambda|_F=0$, that is, $\lambda|_F$ is injective.
 So, $F\stackrel{\lambda|_F}{\cong}L$. It contradicts to
 $dim_kL\stackrel{<}{\not=}dim_kF$.

 From (1)-(8), we finish the description of Example 3.1. Due to
 this example, we know that a general left Artinian algebra with lifted
 quotient cannot be covered by its correspondent $\mathcal A$-path-type tensor algebra.
 This is the reason we introduce pseudo path algebra and $\mathcal A$-path-type pseudo
 tensor algebra so as to replace generalized path algebra and
 $\mathcal A$-path-type tensor algebra in order to cover left Artinian algebras with lifted
 quotients.


 However, there exists still some special class of left Artinian
 algebras which can be covered by its correspondent $\mathcal A$-path-type tensor algebra
 and moreover by a generalized path algebra. This point can be
 seen in the next section, but we will have to restrict a left Artinian
 algebra to be finite dimensional.

\section{When The Radical Is 2-Nilpotent}

In this section, we will always assume the radical $r$ of a finite
dimensional algebra $A$ is $2$-nilpotent, that is, $r^{2}=0$. And,
suppose that $A$ is splitting over its radical $r$ such that
$A=B\oplus r$ with $\widehat{B}=\{B_1,\cdots,B_s\}$ the set of
primitive orthogonal simple subalgebras of $A$ as constructed
above.

For $\bar{x}=x+r\in A/r$, let $\bar{x}\cdot r\stackrel{\rm
def}{=}xr$ and $r\cdot\bar{x}=rx$, then $r$ is a finitely
generated $A/r$-bimodule. For
$A/r=\bigoplus_{i=1}^{s}\overline{A}_{i}$ where $\overline{A}_{i}$
is a simple subalgebra for each $i$, $r$ is a finitely generated
$\overline{A}_{i}$-$\overline{A}_{j}$-bimodule for each pair
$(i,j)$, whose rank is written as $l_{ij}$. Now, $r=A/r\cdot
r\cdot A/r=\sum_{i,j=1}^{s}\overline{A}_{i}\cdot r\cdot
\overline{A}_{j}=\sum_{i,j=1}^{s}$$_{i}M_{j}$ where
$_{i}M_{j}\stackrel{\rm def}{=}\overline{A}_{i}\cdot r\cdot
\overline{A}_{j}$. Then, for $k\not=i$, $\overline{A}_{k}\cdot
_{i}M_{j}=\overline{A}_{k}\cdot(\overline{A}_{i}\cdot r\cdot
\overline{A}_{j})=B_{k}B_{i}rB_{j}=0$, thus $\overline{A}_{k}\cdot
_{i}M_{j}=0$; similarly, for $k\not=j$,
$_{i}M_{j}\cdot\overline{A}_{k}=0$. Hence, we get the
$\mathcal{A}$-path-type tensor algebra $T(A/r,r)$ and the
corresponding $\mathcal{A}$-path algebra
 $k(\Delta,\mathcal{A})$ with $\mathcal{A}=\{\overline{A}_{i}: i\in\Delta_{0}\}$.
 $\Delta$ is called
 the {\em quiver of $A$}.
Similarly with Lemma 3.1, we have the following results:
\begin{lemma}
Assume $A$ is a finite dimensional $k$-algebra with $2$-nilpotent
radical $r=r(A)$ and $A$ is splitting over the radical $r$. Let
$\widehat{B}=\{B_{1},\cdot\cdot\cdot,B_{s}\}$ be
 the set of primitive radical-orthogonal simple subalgebras of $A$ as
constructed above.
 Write $A/r=\bigoplus_{i=1}^{s}\overline{A}_{i}$ where
$\overline{A}_{i}$ are simple algebras for all $i$. The following
statements hold:

 {\em (i)} Let $\{r_{1},\cdot\cdot\cdot,r_{t}\}$ be a set of generators of the $A/r$-bimodule $r$. Then
$B_{1}\cup\cdot\cdot\cdot\cup
B_{s}\cup\{r_{1},\cdot\cdot\cdot,r_{t}\}$ generates $A$ as a
$k$-algebra;

{\em (ii)} There is a surjective $k$-algebra homomorphism
$\widetilde{f}: T(A/r,r)\rightarrow A$ with
\textrm{Ker}$\widetilde{f}=\bigoplus_{j\geq 2}(r)^{j}$, where
$(r)^{j}$ denotes
$r\otimes_{A/r}r\otimes_{A/r}\cdot\cdot\cdot\otimes_{A/r}r$ with
$j$ copies of $r$.
\end{lemma}
{\em Proof}: It is easy to see that $r$ is a $(A/r)$-bimodule with
the action as $\overline{A}_{i}\cdot r=B_{i}r$. Note that
$\overline{A}_{i}\overline{A}_{j}\cdot r=0\cdot r=0$, and on the
other hand, $\overline{A}_{i}\overline{A}_{j}\cdot
r=(B_{i}B_{j}+r)\cdot r=rr=0$, therefore this action is
well-defined. The proof of (i) can be given as similar as the
proof of Lemma 3.1(i) in the case rl$(A)=2$. Next, we prove (ii).

 $r=A/r\cdot r\cdot
A/r=\sum_{i,j=1}^{s}\overline{A}_{i}\cdot
r\cdot\overline{A}_{j}=\sum_{i,j=1}^{s}B_{i}rB_{j}$ is a direct
sum decomposition due to $B_{i}^{2}=B_{i}$ and $B_{i}B_{j}=0$ for
$i\not=j$.

$(A/r)\oplus r$ is a $k$-algebra in which the multiplication is
taken through the $A/r$-module action of $r$ and the
multiplication of $A/r$ and $r$.

 For each pair of integers
$i,j$ with $1\leq i,j\leq s$, choose elements $\{y_{u}^{ij}\}$ a
 $k$-basis in $B_{i}rB_{j}$.  Then
$\bigcup_{i,j=1}^{s}\{y_{u}^{ij}\}$ is a basis for $r$.

Define $f: (A/r)\oplus r\rightarrow A$ by $f|_{r}=id_{r}$, that
is, $f(y_{u}^{ij})=y_{u}^{ij},\;$  and
$f|_{\overline{A}_{i}}=\eta^{-1}$. Then,  $f|_{A/r}$:
$A/r\rightarrow B=f(A/r)$ is a $k$-algebra isomorphism since
$B\stackrel{\eta|_{B}}{\cong}A/r$, and $f|_{r}$: $r\rightarrow
f(r)=r\subset A$ is an embedded homomorphism of
 $A/r$-bimodules.
 Hence by Lemma 2.2, there is a unique
algebra morphism $\widetilde{f}$: $T(A/r,r)\rightarrow A$ such
that $\widetilde{f}|_{(A/r)\oplus r}=f$.

Firstly, $\bigcup_{i,j,u}\{y_{u}^{i,j}\}\subset\widetilde{f}(r)$
and $B_{1}\cup\cdot\cdot\cdot\cup B_{s}\subset\widetilde{f}(A/r)$.
From (i), it follows that $\bigcup_{i,j,u}\{y_{u}^{i,j}\}\cup
B_{1}\cup\cdot\cdot\cdot\cup B_{s}$ generates $A$ as a
$k$-algebra, then $\widetilde{f}$ is surjective. On the other
hand, $f|_{A/r}$ and $f|_{r}$ are monomorphic, so
$\widetilde{f}|_{(A/r)\oplus r}$: $(A/r)\oplus r\rightarrow A$ is
a monomorphism. Then
\textrm{Ker}$\widetilde{f}\subset\bigoplus_{j\geq 2}(r)^{j}$.
 Moreover, $\widetilde{f}((r)^{j})=0$ for $j\geq 2$ since
$r^{j}=0$ in this case. Therefore, $\bigoplus_{j\geq
2}(r)^{j}\subset$\textrm{Ker}$\widetilde{f}$. Thus,
$$\textrm{Ker}\widetilde{f}=\bigoplus_{j\geq 2}(r)^{j}.
\;\;\;\;\;\square$$

In the proof of this lemma, since $f|_r=id_r$, it is naturally an
$A/r$-homomorphism. So, the condition of Lemma 2.2 is satisfied by
$T(A/r,r)$. In general, this is not true for $T(A/r,r)$ in the
case that $r^2\not=0$.

\begin{theorem}
{\em (Generalized Gabriel's Theorem With 2-Nilpotent Radical)}

 Assume that
$A$ is a finite dimensional $k$-algebra with $2$-nilpotent radical
$r=r(A)$ and $A$ is splitting over the radical $r$. Then,
 $A\cong k(\Delta,\mathcal{A},\rho)$ with
$\widetilde J^{2}\subset\langle\rho\rangle\subset\widetilde
J^{2}+\widetilde J\cap$ \textrm{Ker}$\widetilde{\varphi}$ where
$\Delta$ is the quiver of $A$ and $\rho$ is a set of relations of
$k(\Delta,\mathcal{A})$, $\widetilde{\varphi}$ is defined as in
Proposition 2.9.
\end{theorem}
 {\em Proof}:  Let $\Delta$ be the associated quiver of $A$. By Lemma 4.2(ii),
we have the surjective $k$-algebra homomorphism $\widetilde{f}$:
$T(A/r,r)\rightarrow A$. By Proposition 2.9, there is a surjective
$k$-algebra homomorphism $\widetilde{\varphi}$:
$k(\Delta,\mathcal{A})\rightarrow T(A/r,r)$ such that for any
$t\geq 1$, $\widetilde{\varphi}(\widetilde J^{t})=\bigoplus_{j\geq
t}(r)^{j}$. Then $\widetilde{f}\widetilde{\varphi}$:
$k(\Delta,\mathcal{A})\rightarrow A$ is a surjective $k$-algebra
morphism where
$I=$\textrm{Ker}$(\widetilde{f}\widetilde{\varphi})=\widetilde{\varphi}^{-1}(\bigoplus_{j\geq
2}(r)^{j})$ since \textrm{Ker}$\widetilde{f}=\bigoplus_{j\geq
2}(r)^{j}=\widetilde{\varphi}(\widetilde J^{2})$.

As a special case of the correspondent part of the proof of
Theorem 3.2, we have $\widetilde
J^{t}\subset\widetilde{\varphi}^{-1}\widetilde{\varphi}(\widetilde
J^{t})\subset \widetilde J^{t}+\widetilde{\phi}(\bigoplus_{j\leq
t-1}(r)^{Fj})\cap\widetilde{\phi}($\textrm{Ker}$\pi)$ for $t\geq
1$.
  Hence,

$\widetilde
J^{2}\subset\widetilde{\varphi}^{-1}\widetilde{\varphi}(\widetilde
J^{2})=\widetilde{\varphi}^{-1}($\textrm{Ker}$\widetilde{f})=
I\subset\widetilde J^{2}+\widetilde{\phi}(\bigoplus_{j\leq 1}(r)^{Fj})\cap\widetilde{\phi}($\textrm{Ker}$\pi)\subset\widetilde J^{2}+\widetilde J\cap\widetilde{\phi}($\textrm{Ker}$\pi)$.\\
But,
$\widetilde{\phi}($\textrm{Ker}$\pi)=\widetilde{\phi}(\pi^{-1}(0))=\widetilde{\varphi}^{-1}(0)=$\textrm{Ker}$\widetilde{\varphi}$.
Then we get $\widetilde J^{2}\subset I\subset\widetilde
J^{2}+\widetilde J\cap$\textrm{Ker}$\widetilde{\varphi}$.

$\widetilde J^{2}$ is the ideal finitely generated in
$k(\Delta,\mathcal{A})$ by all $\mathcal{A}$-paths of length $2$.
$k(\Delta,\mathcal{A})/\widetilde J^{2}$ is generated finitely as
a $k$-linear space by all $\mathcal{A}$-paths of length less than
$2$. As a $k$-subspace, so is $I/\widetilde J^{2}$, too.
 Then $I$ is a finitely generated ideal
in $k(\Delta,\mathcal{A})$, assume
$\{\sigma_{1},\cdot\cdot\cdot,\sigma_{p}\}$ is its set of finite
generators. Moreover, $\sigma_{l}=\sum_{1\leq i,j\leq
s}\overline{e}_{i}\sigma_{l}\overline{e}_{j}$ where
$\overline{e}_{i}\sigma_{l}\overline{e}_{j}$
 is a relation on the $\mathcal{A}$-path algebra
 $k(\Delta,\mathcal{A})$. Therefore, for
 $\rho=\{\overline{e}_{i}\sigma_{l}\overline{e}_{j}: 1\leq i,j\leq s, 1\leq l\leq p\}$, we get
 $I=\langle\rho\rangle$.
  Hence
  $k(\Delta,\mathcal{A},\rho)=k(\Delta,\mathcal{A})/\langle\rho\rangle\cong A$
  with
  $\langle\rho\rangle=$\textrm{Ker}$(\widetilde{f}\widetilde{\varphi})$ and $\widetilde J^{2}\subset\langle\rho\rangle
  \subset\widetilde J^{2}+\widetilde J\cap$\textrm{Ker}$\widetilde{\varphi}$.
$\;\;\;\;\square$

\begin{corollary}
 Assume that
$A$ is a finite dimensional $k$-algebra with $2$-nilpotent radical
$r=r(A)$. Then,
 $A\cong k(\Delta,\mathcal{A},\rho)$ with
$\widetilde J^{2}\subset\langle\rho\rangle\subset\widetilde J$
where $\Delta$ is the quiver of $A$ and $\rho$ is a set of
relations of $k(\Delta,\mathcal{A})$, if one of the following
conditions holds:

 {\em (i)}$\;\;$ {\em Dim}$A/r\leq 1$ for the radical $r$ of $A$;

 {\em (ii)}$\;\;$ $A/r$ is separable;

{\em (iii)}$\;\;$ $k$ is a perfect field (for example,   when {\em
char}$k=0$ or $k$ is a finite field).
\end{corollary}

As similar as for Theorem 3.5, in the case that
$\langle\rho\rangle\subset J_{\Delta}^2$, there is also the
uniqueness of the correspondent $\mathcal A$-path algebra and
quiver of
 a finite dimensional algebra
  up to isomorphism:
\begin{theorem}
Assume $A$ is a finite dimensional $k$-algebra. Let
$A/r(A)=\bigoplus^{p}_{i=1}\overline{A}_{i}$ with simple algebras
$\overline{A}_{i}$. If there is a quiver $\Delta$ and a
generalized path algebra $k(\Delta,\mathcal{B})$ with a set of
simple algebras $\mathcal{B}=\{B_{1},\cdot\cdot\cdot,B_{q}\}$ and
$\rho$ a set of relations satisfying that $A\cong
k(\Delta,\mathcal{B},\rho)$ with
$J_{\Delta}^{t}\subset\langle\rho\rangle\subset J_{\Delta}^2$ for
some $t$ and $J_{\Delta}$ the ideal in $k(\Delta,\mathcal{B})$
generated by all elements in $k(\Delta_{1},\mathcal{B})$, then
$\Delta$ is the quiver of $A$ and $p=q$ such that
$\overline{A}_{i}\cong B_{i}$ for $i=1,...,p$ after reindexing.
\end{theorem}

This theorem can be proved as similar as for Theorem 3.5, it
suffices to replace $\mathcal A$-path-type tensor algebra and
$\mathcal A$-path with $\mathcal A$-path-type pseudo tensor
algebra and $\mathcal A$-pseudo path respectively.

As said in Section 3, after Theorem 3.5 that in the main results
of Section 3 and 4, $\langle\rho\rangle$ seem not
  admissible. Therefore, the authors think it would be
necessary to
  consider under what condition $\langle\rho\rangle$ is
  admissible such that the conclusions in Theorem 3.5 and 4.4
  can hold.

By Fact 2.5, $\mathcal A$-path-type tensor algebra and $\mathcal
A$-path algebra can be covered respectively by $\mathcal
A$-path-type pseudo tensor algebra and $\mathcal A$-pseudo path
algebra. So, we can also describe the Generalized Gabriel's
Theorem With $2$-Nilpotent Radical through $\mathcal A$-pseudo
path algebra. As a corollary of Theorem 4.3, one has the
following:
\begin{proposition}
 Assume that
$A$ is a finite dimensional $k$-algebra with $2$-nilpotent radical
$r=r(A)$. Then,
 $A\cong PSE_k(\Delta,\mathcal{A},\rho)$ with
$J^{2}\subset\langle\rho\rangle\subset J$ where $\Delta$ is the
quiver of $A$ and $\rho$ is a set of relations on
$PSE_k(\Delta,\mathcal{A})$.
\end{proposition}
{\em Proof:}  We have the composition of surjective homomorphisms:
$PSE_k(\Delta,\mathcal
A)\stackrel{\iota}{\rightarrow}k(\Delta,\mathcal
A)\stackrel{\widetilde{f}\widetilde{\varphi}}{\rightarrow}A$. Then
$A\cong PSE_k(\Delta,\mathcal
A)/\textrm{Ker}(\widetilde{f}\widetilde{\varphi}\iota)$, where
$\textrm{Ker}(\widetilde{f}\widetilde{\varphi}\iota)=\iota^{-1}(\textrm{Ker}(\widetilde{f}\widetilde{\varphi}))$.

By Theorem 4.3, $\widetilde J^{2}\subset
\textrm{Ker}(\widetilde{f}\widetilde{\varphi})\subset\widetilde
J^{2}+\widetilde J\cap$\textrm{Ker}$\widetilde{\varphi}$. Thus,
$$\iota^{-1}(\widetilde J^{2})\subset\iota^{-1}
(\textrm{Ker}(\widetilde{f}\widetilde{\varphi}))\subset\iota^{-1}(\widetilde
J^{2})+\iota^{-1}(\widetilde J\cap
\textrm{Ker}\widetilde{\varphi}).$$ But, since
$\iota^{-1}(\widetilde J)=J$, it follows $\iota^{-1}(\widetilde
J^2)=J^2$ and $\iota^{-1}(\widetilde J\cap
\textrm{Ker}\widetilde{\varphi})\subset J$. Thus, we get
$$J^2\subset \textrm{Ker}(\widetilde{f}\widetilde{\varphi}\iota)\subset
J.$$

By Proposition 2.11 (ii), there is a set $\rho$ of relations on
$PSE_k(\Delta,\mathcal A)$ such that
$\textrm{Ker}(\widetilde{f}\widetilde{\varphi}\iota)=\langle\rho\rangle$.
Then, $A\cong PSE_k(\Delta,\mathcal
A)/\textrm{Ker}(\widetilde{f}\widetilde{\varphi}\iota)=PSE_k(\Delta,\mathcal
A)/\langle\rho\rangle=PSE_k(\Delta,\mathcal A,\rho)$ and
$J^2\subset\langle\rho\rangle\subset J$. $\;\;\;\;\square$
\\

So far, in Section 3 and this section, we have established the
isomorphisms between an algebra and its $\mathcal A$-pseudo path
algebra or $\mathcal A$-path algebra with relations (see Theorem
3.2 and Proposition 4.5) in the cases  either this algebra is left
Artinian with splitting
 over its radical or it is finite-dimensional with $2$-nilpotent.
 However, it seems to be difficult to discuss the same question for an
 arbitrary algebra. Our illusion is whether it would be possible to
 characterize an arbitrary finite-dimensional algebra  through the
 combination of the two methods for a left Artinian algebra with splitting
 over its radical and a finite-dimensional algebra with $2$-nilpotent.

 In fact, for a finite-dimensional algebra $A$, we can start from
 $B=A/r^2$ where $r=r(A)$ the radical of $A$. Consider
 $r(A/r^2)=r/r^2$, denoted as $\widehat{r}$. Then
 $\widehat{r}^2 =r^{2}/r^{2}=0$. By Lemma 4.2(ii), there is a
 surjective homomorphism of algebras $\widetilde{f}:\; T((A/r^2)/(r/r^2),\;r/r^2)\rightarrow
 A/r^2$.

 But, we have $(A/r^2)/(r/r^2)\cong A/r$. So,
 $$\widetilde{f}:\; T(A/r,\;r/r^2)\rightarrow
 A/r^2$$
 is a surjective homomorphism of algebras.

 On the other hand, according to the method in Section 3, in order
 to gain the correspondent Gabriel Theorem for this $A$, the key
 is to give a homomorphism of algebras $\alpha$ as  $\widetilde{f}$ in Lemma
 3.1. Therefore, this question may be thought to find a surjective homomorphism of algebras
 $\alpha$ such that the following diagram commutes:

\[
\begin{diagram}
 \node[2]{T(A/r,\;r/r^2)}\arrow{sw,l}{\alpha}\arrow{s,r}{\widetilde{f}}\\
\node{A}\arrow{e,t}{\pi}\node[1]{A/r^2}\arrow{e,t}{}\node[1]{0}
 \end{diagram}
 \]
 where $\pi$ denotes the natural homomorphism. If such $\alpha$
 exists, the generalized Gabriel Theorem should hold for this
 finite-dimensional algebra $A$.
\\[3.5mm]
{\bf Acknowledgement} {\em This project was supported by the
Program for New Century Excellent Talents in University
(No.04-0522) and the National Natural Science Foundation of China
(No.10571153). The author takes this opportunity to express thanks
to Prof. Shao-xue Liu for his kind advice. His assistance is
important for this work. Also, he would like to
 thank Prof. W. Crawley-Boevey, Prof. Z.Z.Lin, Prof. B.Zhu, Dr. G.X.Liu and Dr. L.L.Chen for their helpful conversations
and suggestions.}
\\
\\
 Fang Li (fangli@zju.edu.cn)
\\
 Department of Mathematics, Zhejiang University, Hangzhou, Zhejiang 310027,  China

\end{document}